\documentclass[aos,seceqn,citesort,dvips]{arximspdf}
\usepackage{dcolumn}
\usepackage{graphicx}

\doi{10.1214/09-AOS785}
\volume{38}
\issue{5}
\pubyear{2010}
\firstpage{3164}
\lastpage{3190}

\makeatletter

\newcolumntype{d}[1]{D{.}{.}{#1}}

\newcommand{\joint}{\,\square\,}

\newtheorem{lemma}{Lemma}[section]
\newtheorem{theorem}{Theorem}[section]
\newtheorem{corollary}[theorem]{Corollary}

\newproclaim{definition}{Definition}

\newproclaim{egg}{Example}

\makeatother

\begin{document}
\begin{frontmatter}

\title{Decomposition tables for experiments. II. Two--one~randomizations}
\runtitle{Decomposition tables II. Two--one randomizations}

\begin{aug}
\author[A]{\fnms{C. J.} \snm{Brien}\corref{}\ead[label=e1]{chris.brien@unisa.edu.au}\ead[label=u1,url]{http://chris.brien.name}} and
\author[B]{\fnms{R. A.} \snm{Bailey}\ead[label=e2]{r.a.bailey@qmul.ac.uk}}
\runauthor{C. J. Brien and R. A. Bailey}
\affiliation{University of South Australia and Queen Mary University
of London}
\address[A]{School of Mathematics and Statistics\\
University of South Australia\\
Mawson Lakes Boulevard\\
Mawson Lakes, SA 5095\\
Australia \\
\printead{e1}\\
\printead{u1}}
\address[B]{School of Mathematical Sciences\\
Queen Mary University of London\\
Mile End Road\\
London E1 4NS\\
United Kingdom\\
\printead{e2}}
\end{aug}

\received{\smonth{7} \syear{2009}}
\revised{\smonth{12} \syear{2009}}

%
\begin{abstract}
We investigate structure for pairs of randomizations that do not
follow each other in a chain. These are unrandomized-inclusive,
independent, coinc\-ident or double randomizations. This involves taking
several structures that satisfy particular relations and combining
them to form the appropriate orthogonal decomposition of the data space
for the experiment.
We show how to establish the decomposition table giving the sources of
variation, their relationships and their degrees of freedom,
so that competing designs can be evaluated. This leads to recommendations
for when the different types of multiple randomization should be used.
\end{abstract}

%
\begin{keyword}[class=AMS]
\kwd[Primary ]{62J10}
\kwd[; secondary ]{62K99}.
\end{keyword}
\begin{keyword}
\kwd{Analysis of variance}
\kwd{balance}
\kwd{decomposition table}
\kwd{design of experiments}
\kwd{efficiency factor}
\kwd{intertier interaction}
\kwd{multiphase experiments}
\kwd{multitiered experiments}
\kwd{orthogonal decomposition}
\kwd{pseudofactor}
\kwd{structure}
\kwd{tier}.
\end{keyword}

\end{frontmatter}

\section{Introduction}
\label{s:introII}

\begin{figure}

\includegraphics{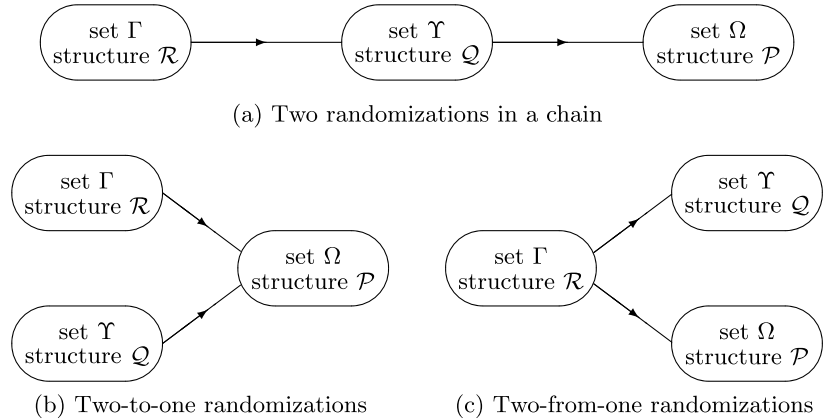}

\caption{The three possibilities for a pair of randomizations.}
\label{fig:intro}
\end{figure}

The purpose of this paper, and its prequel \cite{Brien09b}, is to
establish the orthogonal decomposition of the data space for
experiments that involve multiple randomizations \cite{Brien06},
so that the properties of proposed designs can be evaluated.
In \cite{Brien09b}, this was done for randomizations that follow each
other in a chain, as in Figure \ref{fig:intro}(a).
Here, analogous results to those in \cite{Brien09b} are obtained for
experiments in which the randomizations are two-to-one,
as in Figure \ref{fig:intro}(b).
In such randomizations, two different sets of objects are directly randomized
to a third, as in Figures \ref{fig:cherry}, \ref{fig:indep}
and \ref{fig:plant}.
The unrandomized-inclusive, independent and coincident randomizations
from \cite{Brien06} are of this type. Also covered are experiments in
which the randomization is two-from-one in
that two different sets of objects have a single set of objects randomized
to them; that is, experiments with double randomizations
\cite{Brien06} [see Figure \ref{fig:intro}(c)].

As in \cite{Brien09b}, we always denote the set of observational
units by $\Omega$, so that the data space is the set $V_\Omega$ of all
real vectors indexed by $\Omega$. This data space has an orthogonal
decomposition into subspaces defined by inherent factors and
managerial constraints. We call this decomposition the ``structure''
on $\Omega$, and identify it with the set $\mathcal{P}$ of mutually
orthogonal idempotent matrices which project onto those subspaces.
Thus if $\mathbf{P}\in\mathcal{P}$ then $\mathbf{P}$ is an $\Omega
\times\Omega$ matrix, because its rows and columns are labelled by the
elements of $\Omega$ \cite{Brien09b}.

In the setting of Figure \ref{fig:intro}(b), there are two other sets,
$\Upsilon$ and $\Gamma$, which typically contain treatments of
different types to be randomized to $\Omega$. For example, in
Figure \ref{fig:cherry},
the set of treatments ($\Gamma$) and the set of rootstocks
($\Upsilon$) are randomized to the set of trees ($\Omega$).
Then $V_\Upsilon$ is the space of all real vectors
indexed by $\Upsilon$, and $V_\Gamma$ is defined similarly. Each of
the sets $\Upsilon$ and $\Gamma$ also has a structure
defined on it, the structures being orthogonal decompositions of
$V_\Upsilon$ and $V_\Gamma$, respectively. These are identified with
complete sets $\mathcal{Q}$ and $\mathcal{R}$ of mutually orthogonal
idempotent matrices.

%

There is an immediate technical difficulty. As first defined, a matrix
$\mathbf{Q}$ in $\mathcal{Q}$ is not the same size as a matrix
$\mathbf{P}$ in $\mathcal{P}$. However, the outcome of the
randomization of $\Upsilon$ to $\Omega$ is a function $f$ which allocates
element $f(\omega)$ of $\Upsilon$ to observational unit~$\omega$.
This function defines a subspace $V_\Upsilon^f$ of $V_\Omega$
isomorphic to $V_\Upsilon$. Similarly, the outcome of the
randomization of $\Gamma$ to $\Omega$ is a function $g$ which allocates
element $g(\omega)$ of $\Gamma$ to observational unit $\omega$.
Thus we have a subspace $V_\Gamma^g$ of $V_\Omega$
isomorphic to~$V_\Gamma$. From now on, we identify $V_\Upsilon^f$
with $V_\Upsilon$, and $V_\Gamma^g$ with $V_\Gamma$.
We also assume that equation (4.1) in \cite{Brien09b} holds for both
$f$ and $g$, so that we may regard each matrix $\mathbf{Q}$ in
$\mathcal{Q}$ and each matrix $\mathbf{R}$ in $\mathcal{R}$ as an
$\Omega\times\Omega$ matrix without losing orthogonality or
idempotence. This condition is satisfied for all equi-replicate
allocations, and for many others.

In \cite{Brien09b} it was seen that a standard two-tiered experiment
has just two sets of objects, $\Omega$ and $\Upsilon$ say, typically
observational units and treatments. To evaluate the design for such
an experiment, one needs
the decomposition of the data space $V_\Omega$ that takes into
account both $\mathcal{P}$ and $\mathcal{Q}$.
Brien and Bailey \cite{Brien09b} introduced the notation
$\mathcal{P} \vartriangleright\mathcal{Q}$ for the set of
idempotents for this
decomposition, and established expressions for its elements under the
assumption that $\mathcal{Q}$ is structure balanced in relation
to $\mathcal{P}$.
They exhibited the decomposition in decomposition tables based on
sources corresponding to the elements of $\mathcal{P}$ and $\mathcal{Q}$.

For the idempotents for two sources from different tiers, such as
$\mathbf{P}$ in $\mathcal{P}$ and $\mathbf{Q}$ in $\mathcal{Q}$, we follow
James and Wilkinson \cite{James71} in defining $\mathbf{Q}$ to have
\textit{first-order balance} in relation to $\mathbf{P}$ if there is
a scalar $\lambda_{\mathbf{PQ}}$
such that $\mathbf{QPQ} = \lambda_{\mathbf{PQ}} \mathbf{Q}$.
If this is satisfied and $\lambda_{\mathbf{PQ}} \ne0$, then
$\mathbf{P} \vartriangleright\mathbf{Q}$ is defined in \cite
{Brien09b} to be
$\lambda_{\mathbf{PQ}}^{-1} \mathbf{PQP}$, which is the
matrix of orthogonal projection onto $\operatorname{Im}\mathbf{PQ}$,
the part of the source $\mathbf{P}$ pertaining to the source $\mathbf{Q}$.
The scalar $\lambda_{\mathbf{PQ}}$ is called the \textit{efficiency
factor}; it lies in $[0,1]$ and indicates the proportion of the information
pertaining to the source $\mathbf{Q}$ that is (partially) confounded
with the source $\mathbf{P}$.
Furthermore, a structure $\mathcal{Q}$ is defined in \cite{Brien09b}
to be \textit{structure balanced} in relation to another
structure $\mathcal{P}$ if
(i) all idempotents from $\mathcal{Q}$ have first-order balance in
relation to all idempotents from $\mathcal{P}$; (ii) all pairs of distinct
elements of $\mathcal{Q}$ remain orthogonal when projected onto an
element of
$\mathcal{P}$,
that is, for all $\mathbf{P}$ in $\mathcal{P}$
and all pairs of distinct $\mathbf{Q}_1$ and $\mathbf{Q}_2$ in
$\mathcal{Q}$,
the product \mbox{$\mathbf{Q}_1\mathbf{P}\mathbf{Q}_2 = \mathbf{0}$}.
If $\mathcal{Q}$ is structure
balanced in relation to $\mathcal{P}$, and $\mathbf{P}\in\mathcal
{P}$, then
the residual subspace for $\mathcal{Q}$ in $\operatorname{Im}\mathbf
{P}$ is
just the
orthogonal complement in $\operatorname{Im}\mathbf{P}$ of all the spaces
$\operatorname{Im}\mathbf{PQ}$: its matrix of orthogonal projection
$\mathbf{P}
\vdash\mathcal{Q}$ is given by
%
%
\begin{equation}\label{eq:resid}
\mathbf{P} \vdash\mathcal{Q} = \mathbf{P} -
{\sum_{\mathbf{Q}\in\mathcal{Q}}\hspace*{-2pt}}'
\mathbf{P} \vartriangleright\mathbf{Q},
\end{equation}
where $\sum_{\mathbf{Q}\in\mathcal{Q}}'$ means summation over all
$\mathbf{Q}$ in $\mathcal{Q}$ with $\lambda_{\mathbf{PQ}}\ne0$.

This notation was extended in \cite{Brien09b} to describe the
decomposition for
three-tiered experiments where the two randomizations follow each
other in a chain, as in composed and
randomized-inclusive randomizations \cite{Brien06}
[see Figure \ref{fig:intro}(a)].
This involved
combining the three structures $\mathcal{P}$, $\mathcal{Q}$ and
$\mathcal{R}$ defined on three sets of objects to yield the two
equivalent decompositions $(\mathcal{P} \vartriangleright\mathcal{Q})
\vartriangleright
\mathcal{R}$ and $\mathcal{P} \vartriangleright(\mathcal{Q}
\vartriangleright
\mathcal{R})$. It was seen that the idempotents of these
decompositions could be any of the following forms:
$(\mathbf{P} \vartriangleright\mathbf{Q}) \vartriangleright\mathbf{R}$,
$\mathbf{P} \vartriangleright(\mathbf{Q} \vartriangleright\mathbf{R})$,
$(\mathbf{P} \vartriangleright\mathbf{Q}) \vdash\mathcal{R}$,
$\mathbf{P} \vartriangleright(\mathbf{Q} \vdash\mathcal{R})$,
and
$\mathbf{P} \vdash\mathcal{Q}$,
where $\mathbf{P}$, $\mathbf{Q}$ and $\mathbf{R}$ are idempotents in
$\mathcal{P}$, $\mathcal{Q}$, $\mathcal{R}$, respectively.
In some cases, some idempotents in
$(\mathcal{P} \vartriangleright\mathcal{Q}) \vartriangleright
\mathcal{R}$ may reduce
to idempotents of the form
$\mathbf{P}$, $\mathbf{P} \vartriangleright\mathbf{Q}$, $\mathbf{Q}$,
$\mathbf{R}$ or $\mathbf{Q} \vartriangleright\mathbf{R}$.

In Sections \ref{s:uincl}--\ref{s:onone} of this paper, corresponding
results are obtained for the two-to-one randomizations:
unrandomized-inclusive, independent and coincident randomizations.
It is shown that, in addition to the decompositions above, the
following decompositions occur:
$\mathcal{P} \vartriangleright\mathcal{R}$,
$(\mathcal{P} \vartriangleright\mathcal{R}) \vartriangleright
\mathcal{Q}$
and $(\mathcal{P} \vartriangleright\mathcal{Q}) \joint
(\mathcal{P} \vartriangleright\mathcal{R})$, where ``$\square$'' denotes
``the combination of compatible decompositions'' in a sense defined in
Section \ref{s:onone}.
Also, the list of forms of idempotents is expanded to include:
$\mathbf{P} \vartriangleright\mathbf{R}$,
$\mathbf{P} \vdash\mathcal{R}$,
$(\mathbf{P} \vdash\mathcal{Q}) \vartriangleright\mathbf{R}$ and
$(\mathbf{P} \vdash\mathcal{Q}) \vdash\mathcal{R}$.

Section \ref{s:double} deals with experiments having
the only two-from-one randomization: double randomizations.

There are differences between different
types of multiple randomization
in the reduced forms for the above
idempotents and in the efficiency factors. Section~\ref{s:summary2}
gives recommendations
for when the different types of multiple randomization should be used.
How the results might be applied to
experiments with more than three tiers is outlined in
Section \ref{s:fourII}. We finish in Section \ref{s:discuss} with a
discussion of a number of issues that arise in the decompositions for
multitiered experiments.

\section{Unrandomized-inclusive randomizations}
\label{s:uincl}

In an experiment with
un\-randomized-inclusive randomizations,
$\Upsilon$
is randomized to $\Omega$ in an initial two-tiered experiment. The
unrandomized-inclusive randomization involves a third
set, $\Gamma$, which is randomized to $\Omega$ taking account of the
result of the first randomization. As for
randomized-inclusive
randomizations, the order of the two randomizations is fixed.

\begin{figure}[b]

\includegraphics{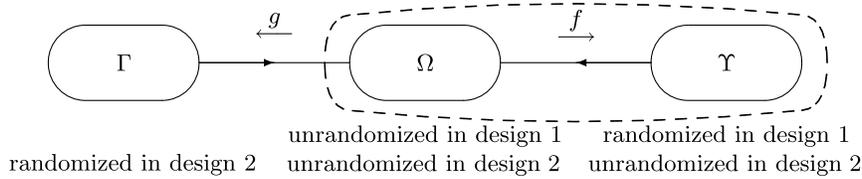}

\caption{Diagram of an experiment with two unrandomized-inclusive
randomizations.}
\label{f:uInclusive}
\end{figure}

Two functions are required to encapsulate the results of these
randomizations, say ${f\colon\Omega\to\Upsilon}$ and ${g\colon
\Omega\to\Gamma}$.
For $\omega$ in $\Omega$, $f(\omega)$ is the element
of $\Upsilon$ assigned to $\omega$ by the first randomization,
and $g(\omega)$ is the element of $\Gamma$ assigned to $\omega$
by the second randomization.
The set-up is represented
diagrammatically in Figure \ref{f:uInclusive}.

%

We consider experiments in which the structure
$\mathcal{Q}$ on $\Upsilon$ is structure balanced in relation to
the structure $\mathcal{P}$ on $\Omega$,
so that the first randomization gives the combined decomposition
$\mathcal{P} \vartriangleright\mathcal{Q}$ of $V_\Omega$
described in \cite{Brien09b}.
The second randomization takes account of $\mathcal{P}
\vartriangleright
\mathcal{Q}$, both in the choice of systematic design and in
restricting the permutations of $\Omega$ to preserve $\mathcal{P}
\vartriangleright\mathcal{Q}$, so we assume that the structure
$\mathcal{R}$ on $\Gamma$
is structure balanced in relation to $\mathcal{P} \vartriangleright
\mathcal{Q}$.

Put $\mathbf{I}_{\mathcal{Q}} =
\sum_{\mathbf{Q}\in\mathcal{Q}}\mathbf{Q}$, which is the matrix of
orthogonal projection onto $V_\Upsilon$. The condition for
$\mathcal{Q}$ to be structure balanced in relation to $\mathcal{P}$ can
be written as $\mathbf{I}_{\mathcal{Q}} \mathbf{P} \mathbf{Q} =
\lambda_{\mathbf{PQ}}\mathbf{Q}$ for all $\mathbf{P}$ in $\mathcal{P}$
and all $\mathbf{Q}$ in $\mathcal{Q}$. Similarly, put
$\mathbf{I}_{\mathcal{R}} = \sum_{\mathbf{R}\in
\mathcal{R}}\mathbf{R}$,
which is the matrix of orthogonal projection onto $V_\Gamma$.
\begin{theorem}\label{th:RbalancedP}
Let $\mathcal{P}$, $\mathcal{Q}$ and $\mathcal{R}$ be orthogonal
decompositions of the spaces $V_\Omega$, $V_{\Upsilon}$ and
$V_{\Gamma}$, respectively, with $V_{\Upsilon} \leq V_\Omega$ and
$V_{\Gamma} \leq V_\Omega$.
If $\mathcal{Q}$ is structure balanced
in relation to $\mathcal{P}$
with efficiency factors $\lambda_{\mathbf{PQ}}$,
and $\mathcal{R}$ is structure balanced
in relation to $\mathcal{P} \vartriangleright\mathcal{Q}$ with
efficiency factors
$\lambda_{\mathbf{P} \vartriangleright\mathbf{Q}, \mathbf{R}}$ and
$\lambda_{\mathbf{P} \vdash\mathcal{Q}, \mathbf{R}}$, then:
\begin{enumerate}[(a)]
\item[(a)]
$\mathcal{R}$ is structure balanced
in relation to $\mathcal{P}$ with
efficiency matrix $\Lambda_{\mathcal{PR}}$ whose entries are
$\lambda_{\mathbf{P} \mathbf{R}} =  ( \lambda_{\mathbf{P}
\vdash
\mathcal{Q}, \mathbf{R}} + \sum_{\mathbf{Q} \in\mathcal{Q}}'
\lambda_{\mathbf{P} \vartriangleright\mathbf{Q}, \mathbf{R}}
)$;\vspace*{1pt}
\item[(b)]
the decomposition
$(\mathcal{P} \vartriangleright\mathcal{Q}) \vartriangleright
\mathcal{R}$ is
\begin{eqnarray*}
& &  \{(\mathbf{P} \vartriangleright\mathbf{Q})
\vartriangleright\mathbf{R}\dvtx\mathbf{P} \in\mathcal{P},  \mathbf
{Q} \in\mathcal{Q},  \mathbf{R} \in\mathcal{R},  \lambda
_{\mathbf{PQ}} \neq0,  \lambda_{\mathbf{P} \vartriangleright
\mathbf{Q}, \mathbf{R}} \neq0 \} \\
&&\qquad{}\cup
 \{(\mathbf{P} \vartriangleright\mathbf{Q}) \vdash\mathcal
{R}\dvtx\mathbf{P} \in\mathcal{P},  \mathbf{Q} \in\mathcal{Q},
\lambda_{\mathbf{PQ}} \neq0 \} \\
&&\qquad{}\cup
 \{(\mathbf{P} \vdash\mathcal{Q}) \vartriangleright\mathbf
{R}\dvtx\mathbf{P} \in\mathcal{P},  \mathbf{R} \in\mathcal{R},
\lambda_{\mathbf{P} \vdash\mathcal{Q}, \mathbf{R}} \neq0 \}
\\
&&\qquad{}\cup
 \{(\mathbf{P} \vdash\mathcal{Q}) \vdash\mathcal{R}\dvtx\mathbf
{P} \in\mathcal{P} \}.
\end{eqnarray*}
\end{enumerate}
\end{theorem}
\begin{pf}
(a)
Because $\mathcal{R}$ is structure balanced in relation to the decomposition
$\mathcal{P} \vartriangleright\mathcal{Q}$, we have
$\mathbf{I}_\mathcal{R} (\mathbf{P} \vartriangleright\mathbf
{Q})\mathbf{R}
= \lambda_{\mathbf{P} \vartriangleright\mathbf{Q}, \mathbf{R}}
\mathbf{R}$
and
$\mathbf{I}_\mathcal{R} (\mathbf{P} \vdash\mathcal{Q}) \mathbf{R}
= \lambda_{\mathbf{P} \vdash\mathcal{Q}, \mathbf{R}} \mathbf{R}$,
for all $\mathbf{P}$ in $\mathcal{P}$, all $\mathbf{Q}$ in $\mathcal{Q}$
with $\lambda_{\mathbf{PQ}}\ne0$, and all $\mathbf{R}$ in $\mathcal{R}$.
Now,
$\mathbf{P} = (\mathbf{P} \vdash\mathcal{Q})
+ \sum_{\mathbf{Q} \in\mathcal{Q}}' \mathbf{P} \vartriangleright
\mathbf{Q}$,
so
\[
\mathbf{I}_\mathcal{R}\mathbf{P}\mathbf{R} =
\mathbf{I}_\mathcal{R} (\mathbf{P} \vdash\mathcal{Q}) \mathbf{R}
+
{\sum_{\mathbf{Q} \in\mathcal{Q}}\hspace*{-2pt}}' \mathbf{I}_\mathcal{R}
(\mathbf{P} \vartriangleright\mathbf{Q}) \mathbf{R}
=  \biggl( \lambda_{\mathbf{P} \vdash\mathcal{Q}, \mathbf{R}}
+ {\sum_{\mathbf{Q} \in\mathcal{Q}}\hspace*{-2pt}}'
\lambda_{\mathbf{P} \vartriangleright\mathbf{Q}, \mathbf{R}} \biggr)
\mathbf{R}.
\]
This proves that $\mathcal{R}$ is structure balanced
in relation to $\mathcal{P}$ with the given efficiency matrix.

(b)
Since $\mathcal{R}$ is structure balanced
in relation to $\mathcal{P} \vartriangleright\mathcal{Q}$, we may
apply the
``$\vartriangleright$'' operator to elements of $\mathcal{P}
\vartriangleright\mathcal{Q}$
and $\mathcal{R}$, to obtain
\[
(\mathbf{P} \vartriangleright\mathbf{Q}) \vartriangleright\mathbf{R}
= \lambda^{-1}_{\mathbf{P} \vartriangleright\mathbf{Q}, \mathbf{R}}
(\mathbf{P} \vartriangleright\mathbf{Q}) \mathbf{R}
(\mathbf{P} \vartriangleright\mathbf{Q})
= \lambda^{-1}_{\mathbf{P} \vartriangleright\mathbf{Q}, \mathbf{R}}
(\lambda^{-1}_{\mathbf{PQ}} \mathbf{PQP}) \mathbf{R}
(\lambda^{-1}_{\mathbf{PQ}} \mathbf{PQP}).
\]
Moreover,\vspace*{1pt} writing $\sum_{\mathbf{R} \in\mathcal{R}}^*$ to mean
summation over $\mathbf{R} \in\mathcal{R}$ with $\lambda_{\mathbf{P}
\vartriangleright\mathbf{Q}, \mathbf{R}} \neq0$,
applying equation (\ref{eq:resid}) to $\mathbf{P} \vartriangleright
\mathbf{Q}$ and $\mathcal{R}$ gives
\[
(\mathbf{P} \vartriangleright\mathbf{Q}) \vdash\mathcal{R}
= \mathbf{P} \vartriangleright\mathbf{Q} - {\sum_{\mathbf{R} \in
\mathcal{R}}\hspace*{-2pt}}^*
(\mathbf{P} \vartriangleright\mathbf{Q}) \vartriangleright\mathbf{R}.
\]
Similarly,
\[
(\mathbf{P} \vdash\mathcal{Q}) \vartriangleright\mathbf{R}
= \lambda^{-1}_{\mathbf{P} \vdash\mathcal{Q}, \mathbf{R}}
(\mathbf{P} \vdash\mathcal{Q}) \mathbf{R} (\mathbf{P} \vdash
\mathcal{Q})
\]
and
\[
(\mathbf{P} \vdash\mathcal{Q}) \vdash\mathcal{R}
= \mathbf{P} \vdash\mathcal{Q} - {\sum_{\mathbf{R} \in\mathcal{R}}\hspace*{-2pt}}^*
(\mathbf{P} \vdash\mathcal{Q}) \vartriangleright\mathbf{R}.
\]
Thus, using Definition 4
in \cite{Brien09b}, the decomposition
$(\mathcal{P} \vartriangleright\mathcal{Q}) \vartriangleright
\mathcal{R}$ is as given.
\end{pf}

The expression for $(\mathcal{P} \vartriangleright\mathcal{Q})
\vartriangleright
\mathcal{R}$ in Theorem \ref{th:RbalancedP}(b) differs from that in
equation (5.1) of \cite{Brien09b}
because
$(\mathbf{P} \vdash\mathcal{Q}) \vartriangleright\mathbf{R}$ is zero
for composed and randomized-inclusive randomizations, but may not be
zero for
unrandomized-inclusive randomizations.

For simplicity, we write the one-dimensional space for the
Mean as $V_0$, with projector $\mathbf{P}_0 = \mathbf{Q}_0 = \mathbf{R}_0
= n^{-1}\mathbf{J}$, where $n =  |\Omega |$ and $\mathbf
{J}$ is the
$n \times n$ all-$1$ matrix.

As Brien and Bailey \cite{Brien06} show, unrandomized-inclusive
randomizations are common in superimposed experiments. In such an
experiment, it may well be the case that $V_\Gamma
\cap V_0^\perp$ is orthogonal to every $\mathbf{P} \vartriangleright
\mathbf{Q}$ of the decomposition $\mathcal{P} \vartriangleright
\mathcal{Q}$. In this case, the decomposition has the simpler form
given by Corollary \ref{eq:u2dec}.
\begin{corollary}
\label{eq:u2dec}
Suppose that $\mathcal{Q}$ is structure balanced in relation to
$\mathcal{P}$ and that $\mathcal{R}$ is structure balanced in relation
to $\mathcal{P} \vartriangleright\mathcal{Q}$.
If $(\mathbf{P} \vartriangleright\mathbf{Q})\mathbf{R} = \mathbf{0}$
for all $\mathbf{P}$ in $\mathcal{P}$, all $\mathbf{Q}$ in
$\mathcal{Q}$ and all $\mathbf{R}$ in $\mathcal{R} \setminus
 \{\mathbf{R}_0 \}$, then
\begin{eqnarray*}
(\mathcal{P} \vartriangleright\mathcal{Q}) \vartriangleright
\mathcal{R}
& =
&  \{(\mathbf{P} \vartriangleright\mathbf{Q})\dvtx\mathbf{P} \in
\mathcal{P},  \mathbf{Q} \in\mathcal{Q},  \lambda_{\mathbf{PQ}}
\neq0 \} \\
& &  {}\cup
 \{(\mathbf{P} \vdash\mathcal{Q}) \vartriangleright\mathbf
{R}\dvtx\mathbf{P} \in\mathcal{P},  \mathbf{R} \in\mathcal{R},
\lambda_{\mathbf{P} \mathbf{R}} \neq0 \}  \\
& &  {}\cup
 \{(\mathbf{P} \vdash\mathcal{Q}) \vdash\mathcal{R}\dvtx\mathbf
{P} \in\mathcal{P} \}.
\end{eqnarray*}
\end{corollary}
\begin{pf}
If $\lambda_{\mathbf{P} \vartriangleright\mathbf{Q}, \mathbf{R}} = 0$
for all $\mathbf{Q} \in\mathcal{Q}$,
then
$\lambda_{\mathbf{P} \mathbf{R}} = \lambda_{\mathbf{P}
\vdash\mathcal{Q}, \mathbf{R}}$.
If this is true for all~$\mathbf{R}$, then
$(\mathbf{P} \vartriangleright\mathbf{Q}) \vdash\mathcal{R} =
\mathbf{P} \vartriangleright\mathbf{Q}$.
The result follows.
\end{pf}
\begin{lemma}
\label{th:qrzero}
Suppose that $\mathcal{Q}$ is structure balanced in relation to
$\mathcal{P}$,
and let $\mathbf{P}\in\mathcal{P} \setminus \{\mathbf
{P}_0 \}$.
The following conditions are equivalent.
\begin{longlist}
\item
$ (\mathbf{P} \vartriangleright\mathbf{Q}) \mathbf{I}_{\mathcal
{R}}={\mathbf
0}$ for all $\mathbf{Q}$ in $\mathcal{Q}$ with
$\lambda_{\mathbf{PQ}} \ne0$.
\item
$\mathbf{QPR} =\mathbf{0}$ for all $\mathbf{Q}$ in $\mathcal{Q}$ and
all $\mathbf{R}$ in $\mathcal{R}$.
\item
$\mathbf{I}_{\mathcal{Q}} \mathbf{P} \mathbf{I}_\mathcal{R} =
\mathbf{0}$.
\end{longlist}
If these are satsified for all $\mathbf{P}$ in $\mathcal{P} \setminus
 \{\mathbf{P}_0 \} $,
then $V_\Upsilon\cap V_0^\perp$
is orthogonal to $V_\Gamma\cap V_0^\perp$,
and all combinations of
elements of $\Upsilon$ with elements of $\Gamma$
occur on $\Omega$.
\end{lemma}
\begin{pf}
If $\lambda_{\mathbf{PQ}} = 0$ then $\mathbf{QP} = \mathbf{0}$ so
$\mathbf{QPR} = \mathbf{0}$. If $\lambda_{\mathbf{PQ}} \ne0$ then
$\mathbf{QPR} = \lambda^{-1}_{\mathbf{PQ}}\times \mathbf{I}_{\mathcal{Q}}
\mathbf{PQP} \mathbf{I}_{\mathcal{R}}\mathbf{R} =
\mathbf{I}_{\mathcal{Q}} (\mathbf{P}\vartriangleright\mathbf{Q})
\mathbf{I}_{\mathcal{R}}\mathbf{R}$. Condition (i) implies that all
these terms are zero, which implies condition (ii).
Summing $\mathbf{QPR}$ over all $\mathbf{Q}$ and all
$\mathbf{R}$ gives
$\mathbf{I}_{\mathcal{Q}} \mathbf{P} \mathbf{I}_{\mathcal{R}}$, so
condition (ii) implies condition (iii). Finally, if
$\lambda_{\mathbf{PQ}} \ne0$ then $(\mathbf{P} \vartriangleright
\mathbf{Q})
\mathbf{I}_\mathcal{R} = \lambda_{\mathbf{PQ}}^{-1} \mathbf{PQP}
\mathbf{I}_\mathcal{R} = \lambda_{\mathbf{PQ}}^{-1} \mathbf{PQ}
(\mathbf{I}_\mathcal{Q} \mathbf{P} \mathbf{I}_\mathcal{R})$,
so condition (iii) implies condition (i).

Summing condition (iii) over all
$\mathbf{P}$ in $\mathcal{P} \setminus \{\mathbf{P}_0 \}$
gives $\mathbf{0} =
\mathbf{I}_\mathcal{Q} (\mathbf{I}_\mathcal{P} - \mathbf{P}_0)
\mathbf{I}_\mathcal{R} = \mathbf{I}_\mathcal{Q} \mathbf
{I}_\mathcal{R}
- \mathbf{I}_\mathcal{Q} \mathbf{P}_0 \mathbf{I}_\mathcal{R} =
(\mathbf{I}_\mathcal{Q} -\mathbf{Q}_0) (\mathbf{I}_\mathcal{R} -
\mathbf{R}_0)$, since $\mathbf{P}_0 = \mathbf{Q}_0 =
\mathbf{R}_0$. This shows that $V_\Upsilon\cap
V_0^\perp$ is orthogonal to $V_\Gamma\cap V_0^\perp$. This implies
that $V_\Upsilon\cap V_\Gamma= V_0$, so
Proposition 2 of
\cite{Bailey96}
shows that the Universe is the only partition marginal to both
$\Upsilon$ and $\Gamma$ considered as factors on $\Omega$. Then
orthogonality and Proposition 3 of \cite{Bailey96}
show that all combinations of $\Upsilon$ and $\Gamma$ occur
on $\Omega$.
\end{pf}

The conditions in Lemma \ref{th:qrzero} are a general form of adjusted
orthogonality \cite{EccRuss75}.
\begin{egg}[(Superimposed experiment in a row-column design)]
\label{eg:SuperStruct}
The initial experiment in Example 10 in \cite{Brien06} is a
randomized complete-block design to investigate cherry rootstocks:
there are three blocks of ten trees each, and there are ten types of
rootstock. Many years later, a set of virus treatments is
superimposed on this, using
the extended Youden square in Table \ref{tab:cherry}.
%
%
\begin{table}[b]
\caption{Extended Youden square showing the
Virus Treatment for each Block--Rootstock combination}\label{tab:cherry}
\begin{tabular*}{\tablewidth}{@{\extracolsep{\fill}}lccccccccccc@{}}
\hline
& & \multicolumn{10}{c@{}}{\textbf{Rootstocks}} \\[-4pt]
& & \multicolumn{10}{c@{}}{\hrulefill}\\
& & \textbf{1} & \textbf{2} & \textbf{3} & \textbf{4}
& \textbf{5} & \textbf{6} & \textbf{7} & \textbf{8} & \textbf{9} & \textbf{10} \\
\hline
Blocks & I & A & B & A & C & D & C & B & E & E & D \\
& II & D & E & B & D & E & A & C & C & A & B \\
& III & E & A & C & E & B & D & D & B & C & A \\
\hline
\end{tabular*}
\end{table}
This ``square'' is a $3\times10$ rectangle whose rows correspond to Blocks
and columns to
Rootstocks. Each of the five treatments occurs twice
in each Block (row), while their disposition in Rootstocks (columns) is
that of a balanced incomplete-block design.
The sets of objects for
this experiment are trees, rootstocks and treatments.
Figure \ref{fig:cherry} shows both randomizations.

%
\begin{figure}

\includegraphics{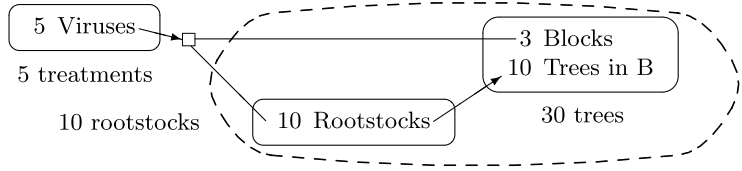}

\caption{Unrandomized-inclusive randomizations in
Example \protect\ref{eg:SuperStruct}: rootstocks are randomized to
trees in the
initial experiment; in the superimposed experiment, treatments are
randomized to trees taking account of the allocation of rootstocks;
$\mathrm{B}$ denotes Blocks.}
\label{fig:cherry}
\end{figure}

For this example,
using the notation for sources in \cite{Brien09b}, but writing
$\mathbf{P}_{\mathrm{Mean}}$ as $\mathbf{P}_0$,
the three structures are
$\mathcal{P} =  \{\mathbf{P}_0, \mathbf{P}_{\mathrm{B}},
\mathbf{P}_{\mathrm{T}[\mathrm{B}]} \}$,
$\mathcal{Q} =  \{\mathbf{Q}_0, \mathbf{Q}_{\mathrm{R}}
\}$ and
$\mathcal{R} =  \{\mathbf{R}_0, \mathbf{R}_{\mathrm{V}}
\}$.
We have
$\mathcal{P} \vartriangleright\mathcal{Q} =
 \{\mathbf{P}_0 \vartriangleright\mathbf{Q}_0, \mathbf
{P}_{\mathrm{B}}, \mathbf{P}_{\mathrm{T}[\mathrm{B}]}
\vartriangleright\mathbf{Q}_{\mathrm{R}}, \mathbf{P}_{\mathrm
{T}[\mathrm{B}]} \vdash\mathcal{Q} \}$,
with
$\mathbf{P}_0 \vartriangleright\mathbf{Q}_0 = \mathbf{P}_0$,
$\mathbf{P}_{\mathrm{T}[\mathrm{B}]} \vartriangleright
\mathbf{Q}_{\mathrm{R}} = \mathbf{Q}_{\mathrm{R}}$ and
$\mathbf{P}_{\mathrm{T}[\mathrm{B}]} \vdash\mathcal{Q} =
\mathbf{P}_{\mathrm{T}[\mathrm{B}]}
- \mathbf{Q}_{\mathrm{R}}$.
See the first two columns of Table~\ref{tab:ANOVASuper}.

%
\begin{table}[b]
\caption{Decomposition table
for Example \protect\ref{eg:SuperStruct}}\label{tab:ANOVASuper}
\begin{tabular*}{\tablewidth}{@{\extracolsep{4in minus 4in}}lrclrcllr@{}}
\hline
\multicolumn{2}{c}{\textbf{trees tier}} &&
\multicolumn{2}{c}{\textbf{rootstocks tier}} &&
\multicolumn{3}{c@{}}{\textbf{treatments tier}} \\[-4pt]
\multicolumn{2}{@{}c}{\hrulefill} && \multicolumn{2}{c}{\hrulefill} &&
\multicolumn{3}{c@{}}{\hrulefill}\\
\multicolumn{1}{@{}l}{\textbf{source}} & \multicolumn{1}{c}{\textbf{d.f.}}
&& \multicolumn{1}{l}{\textbf{source} }&
\multicolumn{1}{l}{\textbf{d.f.}}
&& \multicolumn{1}{l}{\textbf{eff.}} & \textbf{source}
& \multicolumn{1}{l@{}}{\textbf{d.f.}} \\
\hline
Mean & 1 && Mean & 1 & & & Mean & 1\\
\hline
Blocks & 2  \\
\hline
$\mathrm{Trees} [\mathrm{Blocks}]$ & 27  &&
Rootstocks& 9 && $\frac{1}{6}$ & Viruses
& 4 \\
& & & &  &&&Residual & 5\\[-4pt]
& & & \multicolumn{6}{l@{}}{\hrulefill}\\
& && Residual & 18 && $\frac{5}{6}$
&Viruses & 4 \\
& & & & &&& Residual & 14\\
\hline
\end{tabular*}
\end{table}

The efficiency factors for the structure on treatments in relation to
the joint decomposition of trees and rootstocks are
derived from the extended Youden square. Viruses are orthogonal to Blocks,
which means that
$\mathbf{P}_{\mathrm{B}} \mathbf{R}_{\mathrm{V}} = \mathbf{0}$, and
hence $\mathbf{P}_{\mathrm{T[\mathrm{B}]}} \mathbf
{R}_{\mathrm{V}} =
\mathbf{R}_{\mathrm{V}}$.
Viruses have first-order balance in relation
to Rootstocks, with $\lambda_{\mathrm{R}, \mathrm{V}} = 1/6$. Hence
$\mathbf{R}_{\mathrm{V}}(\mathbf{P}_{\mathrm{T}[\mathrm{B}]}
\vartriangleright
\mathbf{Q}_{\mathrm{R}}) \mathbf{R}_{\mathrm{V}} =
\mathbf{R}_{\mathrm{V}}\mathbf{Q}_{\mathrm{R}} \mathbf{R}_{\mathrm
{V}} =
(1/6) \mathbf{R}_{\mathrm{V}}$, and so
$\lambda_{\mathrm{T}[\mathrm{B}]\vartriangleright\mathrm{R},
\mathrm{V}}
= 1/6$. Similarly,
$\lambda_{\mathrm{T}[\mathrm{B}] \vdash
\mathcal{Q}, \mathrm{V}} = 5/6$.
Theorem \ref{th:RbalancedP} shows that the structure on treatments
is orthogonal in relation to the structure on trees since
\[
\lambda_{\mathrm{T} [\mathrm{B}], \mathrm{V}} =
\lambda_{\mathrm{T} [\mathrm{B}]
\vartriangleright\mathrm{R}, \mathrm{V}} +
\lambda_{\mathrm{T}[\mathrm{B}] \vdash
\mathcal{Q}, \mathrm{V}} =
\tfrac{1}{6} + \tfrac{5}{6} = 1.
\]

To obtain the full decomposition
$(\mathcal{P} \vartriangleright\mathcal{Q}) \vartriangleright
\mathcal{R}$,
take $\mathcal{P} \vartriangleright\mathcal{Q}$ and refine it by
$\mathcal{R}$.
In this experiment $V_{\Gamma} \cap V_0^\perp$ is not orthogonal to
$V_\Upsilon\cap V_0^\perp$, because
the Viruses source is not orthogonal to
Rootstocks. This leads to nonorthogonality between
$\mathcal{R}$ and $\mathcal{P} \vartriangleright\mathcal{Q}$.
In particular, the Viruses source
is not orthogonal to $\mathrm{Trees} [\mathrm{Blocks}]
\vartriangleright\mathrm{Rootstocks}$.
Consequently, the decomposition is given by Theorem \ref{th:RbalancedP}(b)
rather than Corollary \ref{eq:u2dec}.
The full decomposition of $V_{\mathrm{trees}}$,
that contains six elements,
one for each line in the decomposition table, is in
Table \ref{tab:ANOVASuper}:
\[
(\mathcal{P} \vartriangleright\mathcal{Q}) \vartriangleright
\mathcal{R} =  \left\{
\begin{array}{l}
(\mathbf{P}_0 \vartriangleright\mathbf{Q}_0)
\vartriangleright\mathbf{R}_0,  \mathbf{P}_{\mathrm{B}},  \\
\bigl(\mathbf{P}_{\mathrm{T}[\mathrm{B}]} \vartriangleright\mathbf
{Q}_{\mathrm{R}}\bigr) \vartriangleright\mathbf{R}_{\mathrm{V}},
\bigl(\mathbf{P}_{\mathrm{T}[\mathrm{B}]} \vartriangleright\mathbf
{Q}_{\mathrm{R}}\bigr) \vdash\mathcal{R},  \\ \bigl(\mathbf{P}_{\mathrm
{T}[\mathrm{B}]} \vdash\mathcal{Q}\bigr) \vartriangleright\mathbf
{R}_{\mathrm{V}},  \bigl(\mathbf{P}_{\mathrm{T}[\mathrm{B}]} \vdash
\mathcal{Q}\bigr) \vdash\mathcal{R}
\end{array}
\right \}
\]
with
\begin{eqnarray*}
(\mathbf{P}_0 \vartriangleright\mathbf{Q}_0) \vartriangleright
\mathbf{R}_0
&=& \mathbf{P}_0 = \mathbf{Q}_0 = \mathbf{R}_0, \\
\bigl(\mathbf{P}_{\mathrm{T}[\mathrm{B}]} \vartriangleright\mathbf
{Q}_{\mathrm{R}}\bigr) \vartriangleright\mathbf{R}_{\mathrm{V}}
&=& 6\bigl(\mathbf{P}_{\mathrm{T}[\mathrm{B}]}\vartriangleright\mathbf
{Q}_{\mathrm{R}}\bigr)\mathbf{R}_{\mathrm{V}}
\bigl(\mathbf{P}_{\mathrm{T}[\mathrm{B}]}\vartriangleright\mathbf
{Q}_{\mathrm{R}}\bigr) \\
&=&
6\mathbf{Q}_{\mathrm{R}} \mathbf{R}_{\mathrm{V}} \mathbf
{Q}_{\mathrm{R}}
= \mathbf{Q}_{\mathrm{R}} \vartriangleright\mathbf{R}_{\mathrm
{V}}, \\
\bigl(\mathbf{P}_{\mathrm{T}[\mathrm{B}]} \vartriangleright\mathbf
{Q}_{\mathrm{R}}\bigr) \vdash\mathcal{R}
&=& \mathbf{P}_{\mathrm{T}[\mathrm{B}]} \vartriangleright\mathbf
{Q}_{\mathrm{R}}
- \bigl(\mathbf{P}_{\mathrm{T}[\mathrm{B}]} \vartriangleright\mathbf
{Q}_{\mathrm{R}}\bigr) \vartriangleright\mathbf{R}_{\mathrm{V}}
\\
&=& \mathbf{Q}_{\mathrm{R}} - \mathbf{Q}_{\mathrm{R}}
\vartriangleright
\mathbf{R}_{\mathrm{V}},\\
\bigl(\mathbf{P}_{\mathrm{T}[\mathrm{B}]} \vdash\mathcal{Q}\bigr)
\vartriangleright\mathbf{R}_{\mathrm{V}}
&=& \lambda^{-1}_{\mathrm{T}[\mathrm{B}]\vdash\mathcal
{Q},\mathrm{V}}
\bigl( \mathbf{P}_{\mathrm{T}[\mathrm{B}]} \vdash\mathcal
{Q}\bigr)\mathbf{R}_{\mathrm{V}}
\bigl(\mathbf{P}_{\mathrm{T}[\mathrm{B}]} \vdash\mathcal{Q}\bigr),
\\
\bigl(\mathbf{P}_{\mathrm{T}[\mathrm{B}]} \vdash\mathcal{Q}\bigr)
\vdash\mathcal{R}
&=& \mathbf{P}_{\mathrm{T}[\mathrm{B}]} \vdash\mathcal{Q}
- \bigl(\mathbf{P}_{\mathrm{T}[\mathrm{B}]} \vdash\mathcal{Q}\bigr)
\vartriangleright\mathbf{R}_{\mathrm{V}}.
\end{eqnarray*}

As expected, this decomposition does contain a nontrivial
idempotent of the form $(\mathbf{P} \vartriangleright\mathbf{Q})
\vartriangleright
\mathbf{R}$.
Also, unlike the chain randomizations in \cite{Brien09b},
it contains an idempotent of the form
$(\mathbf{P} \vdash\mathcal{Q}) \vartriangleright\mathbf{R}$.

The efficiency factors are recorded
in the decomposition
in Table \ref{tab:ANOVASuper},
which shows that the Viruses source
is partly confounded with
both Rootstocks and the part of $\mathrm{Trees}[\mathrm{Blocks}]$ that is
orthogonal to Rootstocks. A consequence of this is that four Rootstocks degrees
of freedom cannot be separated from Virus differences. However, there
are five
Rootstocks degrees of freedom that are orthogonal to Virus differences.
Further, while the Viruses source has first-order balance
in relation to Rootstocks, the reverse is not true.
\end{egg}

\section{Independent or coincident randomizations}
\label{s:onone}

\begin{figure}

\includegraphics{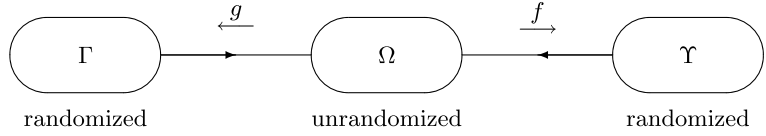}

\caption{Diagram of an experiment with two independent
or two coincident randomizations.}
\label{f:Coincide}
\end{figure}

For independent or coincident randomizations, two sets of objects are
randomized to the third; thus
we could have $\Gamma$ and $\Upsilon$ randomized to $\Omega$. Two
functions are needed to encapsulate the results of these
randomizations, say ${f\colon\Omega\to\Gamma}$ and ${g\colon
\Omega\to\Upsilon}$.
The set-up is represented diagrammatically
in Figure \ref{f:Coincide}.
A particular feature of these randomizations
is that there is no intrinsic ordering of $\Gamma$ and $\Upsilon$,
because neither randomization takes account of the outcome of the other.
Associated with $\Omega$, $\Upsilon$ and $\Gamma$ are the
decompositions $\mathcal{P}$, $\mathcal{Q}$ and $\mathcal{R}$.
We assume that $\mathcal{Q}$ and $\mathcal{R}$ are both structure
balanced in relation to $\mathcal{P}$.

%

%
\begin{figure}[b]

\includegraphics{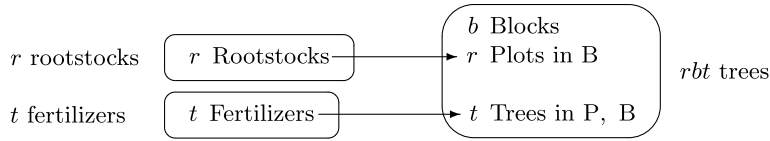}

\caption{Independent randomizations in Example \protect\ref{eg:SuperSplit}:
rootstocks are randomized to trees in such a way that all trees in
each plot have a single type of rootstock; later, fertilizers are
randomized to trees in such a way that each fertilizer is applied to one
tree per plot; $\mathrm{B}$, $\mathrm{P}$ denote Blocks, Plots,
respectively.}
\label{fig:indep}
\end{figure}

The difference between coincident and independent randomizations is
that, for coincident randomizations, there are sources from the two
randomized tiers which are both (partly) confounded with the same
source in the unrandomized tier. For independent randomizations this
does not occur (apart from the Mean).

\subsection{Independent randomizations}
\label{s:indep}
For a pair of independent randomizations,
the two functions are randomized by two permutations
chosen independently from the same group of
permutations of $\Omega$.
The precise definition of independence,
which we were unable to give in \cite{Brien06}, is that
the conditions in Lemma \ref{th:qrzero} are satisfied,
for all $\mathbf{P}$ in $\mathcal{P} \setminus \{\mathbf
{P}_0 \}$,
for all possible outcomes of the two randomizations.
If $\lambda_{\mathbf{PQ}}\lambda_{\mathbf{PR}} \ne0$ then
some outcomes will have $\mathbf{QPR}\ne\mathbf{0}$,
violating these conditions.
Hence independent randomizations require
that $\lambda_{\mathbf{PQ}}\lambda_{\mathbf{PR}} = 0$ for all
$\mathbf{Q}$ in $\mathcal{Q}$ and all $\mathbf{R}$ in $\mathcal{R}$
unless $\mathbf{P}= \mathbf{P}_0$.
Lemma \ref{th:qrzero} shows that,
if $\mathcal{Q}$ and $\mathcal{R}$ are both
structure balanced in relation to $\mathcal{P}$,
then they are also structure balanced
in relation to $\mathcal{P} \vartriangleright\mathcal{R}$
and $\mathcal{P} \vartriangleright\mathcal{Q}$, respectively, with
$\lambda_{\mathbf{P} \vartriangleright\mathbf{Q}, \mathbf{R}} =
\lambda_{\mathbf{P} \vartriangleright\mathbf{R}, \mathbf{Q}} = 0$
unless $\mathbf{P} = \mathbf{P}_0$, $\mathbf{Q} = \mathbf{Q}_0$ and
$\mathbf{R} = \mathbf{R}_0$. Therefore
%
%
\begin{eqnarray}\label{eq:decindep}
(\mathcal{P} \vartriangleright\mathcal{Q}) \vartriangleright
\mathcal{R} &=&
(\mathcal{P} \vartriangleright\mathcal{R}) \vartriangleright
\mathcal{Q}\nonumber\\
&=&  \{\mathbf{P} \vartriangleright\mathbf{Q}\dvtx\mathbf{P}\in
\mathcal{P}, \mathbf{Q}\in\mathcal{Q},  \lambda_{\mathbf
{PQ}}\ne0 \} \nonumber\\
&&  {} \cup \{\mathbf{P} \vdash\mathcal{Q}\dvtx\mathbf
{P} \in\mathcal{P}, \mathbf{P}\mathbf{I}_\mathcal{Q} \ne\mathbf
{0} \}
\nonumber\\[-9pt]\\[-9pt]
&&  {} \cup \{\mathbf{P} \vartriangleright\mathbf
{R}\dvtx\mathbf{P}\in\mathcal{P}, \mathbf{R}\in\mathcal{R},
\lambda_{\mathbf{PR}}\ne0 \} \nonumber\\
&&  {} \cup \{\mathbf{P} \vdash\mathcal{R}\dvtx\mathbf
{P} \in\mathcal{P}, \mathbf{P}\mathbf{I}_\mathcal{R} \ne\mathbf
{0} \}
\nonumber\\
&&  {} \cup \{\mathbf{P}\dvtx\mathbf{P} \in\mathcal
{P}, \mathbf{P}\mathbf{I}_\mathcal{Q} = \mathbf{P}\mathbf
{I}_\mathcal{R} = \mathbf{0} \}.\nonumber
\end{eqnarray}

As outlined in \cite{Brien06}, Section 8.5,
wherever possible we reduce two independent randomizations to a
single randomization. However, as noted in
\cite{Brien06}, Section 4.3,
this is not always possible---for
example, when it is not physically possible to do them simultaneously.
\begin{egg}[(Superimposed experiment using split plots)]
\label{eg:SuperSplit}
Example 6 in \cite{Brien06} is a superimposed experiment
in which the second set of treatments (fertilizers) is
randomized to subunits
(trees) of the original experimental units (plots). The randomizations are
independent, being carried out at different times and with the later
one taking
no account of the earlier one except to force fertilizers to be
orthogonal to
rootstocks. See Figure \ref{fig:indep}. Table \ref{tab:indep} shows the
decomposition.

%
\begin{table}
\caption{Decomposition table
for Example \protect\ref{eg:SuperSplit}}\label{tab:indep}
\begin{tabular*}{\tablewidth}{@{\extracolsep{4in minus 4in}}lclclc@{}}
\hline
\multicolumn{2}{@{}c}{\textbf{trees tier}} & \multicolumn
{2}{c}{\textbf{rootstocks tier}} &
\multicolumn{2}{c@{}}{\textbf{fertilizers tier}} \\[-4pt]
\multicolumn{2}{@{}l}{\hrulefill} & \multicolumn
{2}{c}{\hrulefill} & \multicolumn{2}{c@{}}{\hrulefill} \\
\textbf{source} & \textbf{d.f.} & \textbf{source} & \textbf{d.f.}
& \textbf{source} & \textbf{d.f.} \\
\hline
Mean & 1 & Mean & 1 & Mean & 1\\
\hline
Blocks & $b-1$ \\
\hline
$\mathrm{Plots}[\mathrm{B}]$ &$b(r-1)$ &
Rootstocks & $r-1$  \\
& &  Residual & $(b-1)(r-1)$ \\
\hline
$\mathrm{Trees}[\mathrm{P}\wedge\mathrm{B}]$ &
$br(t-1)$
& & &  Fertilizers & $t-1$ \\
& & & & Residual & $(br-1)(t-1)$\\
\hline
\end{tabular*}
\end{table}

In this example
the independence of the randomizations implies that
$(\mathbf{P}_{\mathrm{P}[\mathrm{B}]} \vartriangleright
\mathbf{Q}_{\mathrm{R}}) \vartriangleright\mathbf{R}_{\mathrm{F}} =
\mathbf{0}$ and so
the conditions in Lemma \ref{th:qrzero} are satisfied.
\end{egg}

\subsection{Coincident randomizations}
\label{s:coincide}

For coincident randomizations, there are idempotents
$\mathbf{P}$ in $\mathcal{P} \setminus \{\mathbf{P}_0 \}$,
$\mathbf{Q}$ in $\mathcal{Q} \setminus \{\mathbf{Q}_0 \}$
and $\mathbf{R}$ in $\mathcal{R} \setminus \{\mathbf
{R}_0 \}$
such that $\mathbf{PQ}$ and $\mathbf{PR}$ are both nonzero. If
$\operatorname{Im}
\mathbf{PQ}$ and $\operatorname{Im}\mathbf{PR}$ are both proper
subspaces of
$\operatorname{Im}\mathbf{P}$, then the relationship between $\mathbf
{Q}$ and
$\mathbf{R}$ depends on the choice of the two independent permutations
used in randomizing $\Upsilon$ and $\Gamma$ to $\Omega$; restricting
one of the randomizations to preserve the relationship would make the
multiple randomizations unrandomized inclusive rather than
coincident. On the other hand, if $\operatorname{Im}\mathbf{PQ} =
\operatorname{Im}
\mathbf{P}$ then $\operatorname{Im}\mathbf{PR}$ is always contained
in $\operatorname{Im}
\mathbf{PQ}$. If $\mathcal{Q}$ is structure balanced in relation
to $\mathcal{P}$ and $\operatorname{Im}\mathbf{PQ} = \operatorname
{Im}\mathbf{P}$,
then $\mathbf{P}\vartriangleright\mathbf{Q} = \mathbf{P}$ and
the two sources corresponding to $\mathbf{Q}$ and $\mathbf{P}$
have the same number of degrees of freedom.
The condition for coincident randomizations hinted at in
\cite{Brien06}, Section 4.2, is precisely that
%
%
\begin{equation}\label{eq:coincond}
\begin{tabular}{p{325pt}}
for all $\mathbf{P}$ in $\mathcal{P}$, $\mathbf{Q}$ in $\mathcal{Q}$
and $\mathbf{R}$ in $\mathcal{R}$,
if $\mathbf{PQ}$ and
$\mathbf{PR}$ are both nonzero then one of $\mathbf{P}
\vartriangleright
\mathbf{Q}$ and $\mathbf{P} \vartriangleright\mathbf{R}$ is equal
to $\mathbf{P}$.
\end{tabular}\hspace*{-35pt}
\end{equation}
A special, commonly occurring, case arises when
$\mathcal{Q}$ and $\mathcal{R}$ can be assigned to the two randomized
sets of objects such that
the following condition is satisfied:
%
%
\begin{equation}\label{eq:coin}
\begin{tabular}{p{310pt}}
for all $\mathbf{P}$ in $\mathcal{P}$ and $\mathbf{Q}$ in $\mathcal{Q}$,
if $\mathbf{PQ}$ and $\mathbf{P}\mathbf{I}_{\mathcal{R}}$ are both
nonzero then\break $\mathbf{P}\vartriangleright\mathbf{Q} = \mathbf{P}$.
\end{tabular}\hspace*{-30pt}
\end{equation}
\begin{theorem}
\label{th:RbalancedPQ}
If $\mathcal{Q}$ and $\mathcal{R}$ are both structure balanced
in relation to $\mathcal{P}$ and condition (\ref{eq:coin}) is satisfied
then $\mathcal{R}$ is structure balanced
in relation to
$\mathcal{P} \vartriangleright\mathcal{Q}$, with $\lambda_{\mathbf
{P}\vartriangleright
\mathbf{Q}, \mathbf{R}} = \lambda_{\mathbf{PR}}$ if
$\lambda_{\mathbf{PQ}} \ne0$ and
$\lambda_{\mathbf{P}\vdash\mathcal{Q}, \mathbf{R}}=
\lambda_{\mathbf{PR}}$ if $\mathbf{P}\vdash\mathcal{Q} \ne
\mathbf{0}$.
Moreover, the decomposition $(\mathcal{P}
\vartriangleright\mathcal{Q}) \vartriangleright\mathcal{R}$ is
%
%
\begin{eqnarray} \label{eq:cdec}
& &
 \{\mathbf{P} \vartriangleright\mathbf{R}\dvtx\mathbf{P} \in
\mathcal{P},  \mathbf{R} \in\mathcal{R},  \lambda_{\mathbf{PR}}
\neq0 \}\nonumber\\
&&\qquad{}\cup
 \{(\mathbf{P} \vdash\mathcal{R})\dvtx\mathbf{P} \in\mathcal
{P},  \mathbf{PI}_{\mathcal{R}}\ne0 \} \nonumber\\[-8pt]\\[-8pt]
&&\qquad{}\cup
 \{(\mathbf{P} \vartriangleright\mathbf{Q}) \dvtx\mathbf{P} \in
\mathcal{P},  \mathbf{Q} \in\mathcal{Q},  \lambda_{\mathbf{PQ}}
\neq0,  \mathbf{PI}_{\mathcal{R}}=0 \}
\nonumber\\
&&\qquad{} \cup
 \{(\mathbf{P} \vdash\mathcal{Q})\dvtx\mathbf{P} \in\mathcal
{P},  \mathbf{PI}_{\mathcal{R}}=0 \}.
\nonumber
\end{eqnarray}
\end{theorem}
\begin{pf}
If $\lambda_{\mathbf{PR}} =0$ then $\mathbf{PR} = \mathbf{0}$,
so $(\mathbf{P} \vartriangleright\mathbf{Q}) \mathbf{R} = \mathbf
{0}$ for all
$\mathbf{Q}$ with $\lambda_{\mathbf{PQ}} \ne0$, and hence
$(\mathbf{P}\vdash\mathcal{Q}) \mathbf{R} = \mathbf{0}$.
Suppose that $\mathbf{P} \mathbf{I}_\mathcal{R} \ne\mathbf{0}$.
Then either $\mathbf{P} \mathbf{I}_{\mathcal{Q}} =
\mathbf{0}$ or there is a unique $\mathbf{Q}$ in $\mathcal{Q}$ with
$\lambda_{\mathbf{PQ}} \ne0$, which satisfies $\mathbf{P} =
\mathbf{P} \vartriangleright\mathbf{Q}$. In the first case,
$\mathbf{P} = \mathbf{P} \vdash\mathcal{Q}$: therefore
$\mathbf{I}_\mathcal{R} (\mathbf{P} \vdash\mathcal{Q}) \mathbf{R} =
\mathbf{I}_\mathcal{R} \mathbf{P} \mathbf{R} = \lambda_{\mathbf{PR}}
\mathbf{R}$,
and $(\mathbf{P}\vdash\mathcal{Q}) \vartriangleright\mathbf{R} =
\mathbf{P} \vartriangleright\mathbf{R}$. In the second case,
$\mathbf{I}_\mathcal{R} (\mathbf{P} \vartriangleright\mathbf{Q})
\mathbf
{R} =
\mathbf{I}_\mathcal{R} \mathbf{P} \mathbf{R} =
\lambda_{\mathbf{PR}} \mathbf{R}$
and $(\mathbf{P}\vartriangleright\mathbf{Q}) \vartriangleright
\mathbf{R} =
\mathbf{P} \vartriangleright\mathbf{R}$.
\end{pf}
\begin{egg}[(A plant experiment)]
\label{eg:PlantExp}
Example 5 in \cite{Brien06} is an experiment to investigate
five varieties and two spray regimes.
Each bench has one spray regime and two seedlings of each variety.
See Figure \ref{fig:plant}. The sets are
positions, seedlings and regimes.
The diagram includes the pseudofactor $\mathrm{S}_1$ for
{$\mathrm{Seedlings}[\mathrm{Varieties}]$},
which indexes the groups of
seedlings randomized to the different benches. Although
the factor
Seedlings is nested in Varieties, $\mathrm{S}_1$ is not, because
each of its levels is taken across all levels of Varieties.

%
%
\begin{figure}[b]

\includegraphics{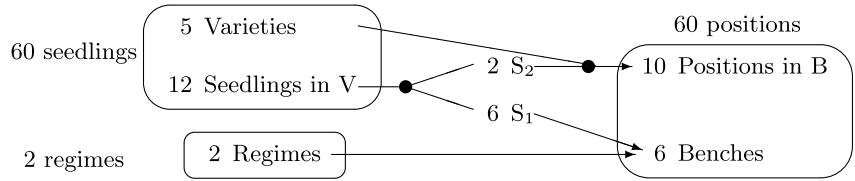}

\caption{Coincident randomizations in Example \protect\ref{eg:PlantExp}:
seedlings and regimes are both randomized to positions; $\mathrm{V}$
denotes Varieties, $\mathrm{B}$ denotes Benches; $\mathrm{S}_1$ and
$\mathrm{S}_2$ are pseudofactors for Seedlings.}
\label{fig:plant}
\end{figure}

The Hasse diagrams displaying the structures for this experiment are
in Figure~\ref{f:HassePlant}. The decomposition
is in Table \ref{tab:ANOVAPlant},
where the source $\mathrm{Seedlings} [\mathrm{Varieties}]
\vdash{\mathrm{S}_1}$ is the part of $\mathrm{Seedlings} [\mathrm
{Varieties}]$ which is orthogonal to the source $\mathrm{S}_1$.

%
\begin{figure}

\includegraphics{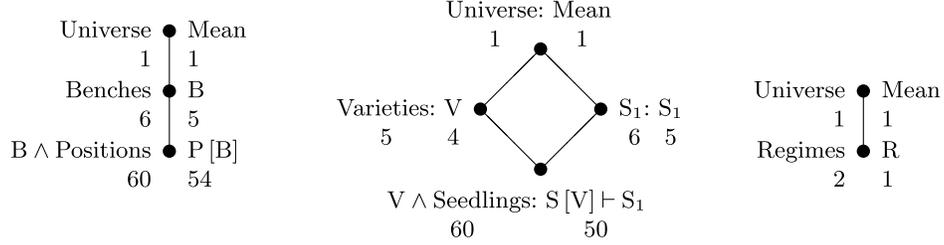}

\caption{Hasse diagrams for Example \protect\ref{eg:PlantExp}.}
\label{f:HassePlant}
\end{figure}

%
\begin{table}[b]
\caption{Decomposition table for
Example \protect\ref{eg:PlantExp}}\label{tab:ANOVAPlant}
\begin{tabular*}{\tablewidth}{@{\extracolsep{4in minus 4in}}lrlrlr@{}}
\hline
\multicolumn{2}{@{}c}{\textbf{positions tier}} &
\multicolumn{2}{c}{\textbf{seedlings tier}} &
\multicolumn{2}{c@{}}{\textbf{regimes tier}}
\\[-4pt]
\multicolumn{2}{@{}l}{\hrulefill} & \multicolumn
{2}{c}{\hrulefill} & \multicolumn{2}{c@{}}{\hrulefill} \\
\textbf{source} & \multicolumn{1}{l}{\textbf{d.f.}}
& \multicolumn{1}{l}{\textbf{source}} &
\multicolumn{1}{l}{\textbf{d.f.}}
& \multicolumn{1}{l}{\textbf{source}} & \multicolumn{1}{l@{}}{\textbf{d.f.}} \\
\hline
Mean & \phantom{0}1 & Mean & \phantom{0}1 & Mean & 1\\
\hline
Benches & \phantom{0}5 & $\mathrm{S}_1$ & \phantom{0}5 & Regimes &1 \\
& & & & Residual &4\\
\hline
$\mathrm{Positions}[\mathrm{Benches}]$ & 54  &
Varieties & \phantom{0}4
\\
& & $\mathrm{Seedlings} [\mathrm{Varieties}]
\vdash\mathrm{S}_1$
& 50\\
\hline
\end{tabular*}
\end{table}

The full decomposition of $V_{\mathrm{positions}}$ in this case contains
five elements and is
\[
(\mathcal{P} \vartriangleright\mathcal{Q}) \vartriangleright
\mathcal{R} =  \left\{
\begin{array}{l}
(\mathbf{P}_0 \vartriangleright\mathbf{Q}_0)
\vartriangleright\mathbf{R}_0,  \\ (\mathbf{P}_{\mathrm{B}}
\vartriangleright\mathbf{Q}_{\mathrm{S}_1}) \vartriangleright
\mathbf{R}_{\mathrm{R}},  (\mathbf{P}_{\mathrm{B}}
\vartriangleright\mathbf{Q}_{\mathrm{S}_1}) \vdash\mathcal{R},
\\ \mathbf{P}_{\mathrm{P}[\mathrm{B}]} \vartriangleright\mathbf
{Q}_{\mathrm{V}},  \mathbf{P}_{\mathrm{P}[\mathrm{B}]}
\vartriangleright\mathbf{Q}_{\mathrm{S}[\mathrm{V}] \vdash\mathrm
{S}_1}
\end{array}
\right\}
.
\]

This experiment clearly meets condition (\ref{eq:coin}),
because the only source for positions which
is nonorthogonal to sources from both of the randomized tiers is
the Benches source,
and the five-dimensional pseudosource $\mathrm{S}_1$ is equal
to Benches.
That is, $\mathbf{P}_{\mathrm{B}} \vartriangleright\mathbf
{Q}_{\mathrm
{S}_1} = \mathbf{P}_{\mathrm{B}}
= \mathbf{Q}_{\mathrm{S}_1}$.
The other source nonorthogonal to Benches is the
one-dimensional source
Regimes, which is a proper subspace of the Benches source,
and so $(\mathbf{P}_{\mathrm{B}} \vartriangleright\mathbf
{Q}_{\mathrm
{S}_1}) \vartriangleright\mathbf{R}_{\mathrm{R}}
= \mathbf{P}_{\mathrm{B}} \vartriangleright\mathbf{R}_{\mathrm{R}} =
\mathbf{R}_{\mathrm{R}}$.

Consequently, the elements of the full decomposition can be written as follows:
\[
(\mathcal{P} \vartriangleright\mathcal{Q}) \vartriangleright
\mathcal{R} =  \left\{
\begin{array}{l}
\mathbf{P}_0,  \mathbf{P}_{\mathrm{B}}
\vartriangleright\mathbf{R}_{\mathrm{R}}, \mathbf{P}_{\mathrm
{B}} \vdash\mathcal{R},  \\ \mathbf{P}_{\mathrm{P}[\mathrm{B}]}
\vartriangleright\mathbf{Q}_{\mathrm{V}},  \mathbf{P}_{\mathrm
{P}[\mathrm{B}]} \vartriangleright\mathbf{Q}_{\mathrm{S}[\mathrm
{V}] \vdash\mathrm{S}_1}
\end{array}
\right\}
.
\]

On noting that $\mathbf{P}_{\mathrm{B}}\mathbf{R}_{\mathrm{R}}
=\mathbf{R}_{\mathrm{R}}$,
$\mathbf{P}_{\mathrm{P}[\mathrm{B}]}\mathbf{Q}_{\mathrm{V}}
= \mathbf{Q}_{\mathrm{V}}$ and
$\mathbf{P}_{\mathrm{P}[\mathrm{B}]}\mathbf{Q}_{\mathrm
{S}[\mathrm{V}]
\vdash\mathrm{S}_1} = \mathbf{Q}_{\mathrm{S}[\mathrm{V}]
\vdash\mathrm{S}_1}$, the decomposition further reduces to
\[
(\mathcal{P} \vartriangleright\mathcal{Q}) \vartriangleright
\mathcal{R}
=  \bigl\{\mathbf{P}_0,  \mathbf{R}_{\mathrm{R}},  \mathbf
{P}_{\mathrm{B}} - \mathbf{R}_{\mathrm{R}},  \mathbf{Q}_{\mathrm
{V}},  \mathbf{Q}_{\mathrm{S}[\mathrm{V}] \vdash\mathrm
{S}_1} \bigr\}.
\]
\end{egg}

Decomposition (\ref{eq:cdec}) is convenient for algorithms, because it is
$(\mathcal{P} \vartriangleright\mathcal{Q}) \vartriangleright
\mathcal{R}$, like
the decompositions in Theorem 5.1(d)
in \cite{Brien09b}, Theorem \ref{th:RbalancedP}(b),
Corollary \ref{eq:u2dec}
and equation (\ref{eq:decindep}). However,
it gives the false impression that the decomposition of $V_\Omega$
must have $\mathcal{P}$ refined by $\mathcal{Q}$, then
$\mathcal{P} \vartriangleright\mathcal{Q}$ refined by $\mathcal{R}$,
suggesting that $\mathcal{Q}$ and $\mathcal{R}$ have different roles.
Moreover, condition (\ref{eq:coin}) does not hold for all pairs of
coincident randomizations.
We therefore introduce another joint
decomposition that emphasizes the symmetry between
$\mathcal{Q}$ and $\mathcal{R}$.
\begin{definition}
Let $\mathcal{B}$ and $\mathcal{C}$ be orthogonal decompositions of
the same space~$V_\Omega$. Then $\mathcal{B}$ is \textit{compatible}
with $\mathcal{C}$ if $\mathbf{BC} = \mathbf{CB}$
for all $\mathbf{B}$ in $\mathcal{B}$ and all $\mathbf{C}$
in $\mathcal{C}$.
\end{definition}

Lemma 2.4 in \cite{Bailey04a}
shows that if $\mathcal{B}$ and
$\mathcal{C}$ are compatible then the nonzero products $\mathbf{BC}$,
for $\mathbf{B}$ in $\mathcal{B}$ and $\mathbf{C}$ in $\mathcal{C}$,
give another orthogonal decomposition of $V_\Omega$, which is a
refinement of both $\mathcal{B}$ and $\mathcal{C}$.
\begin{definition}
If $\mathcal{B}$ and $\mathcal{C}$ are orthogonal decompositions
of $V_\Omega$ which are compatible with each other, then the
decomposition $\mathcal{B} \joint\mathcal{C}$
of $V_\Omega$ is defined to be
\[
\mathcal{B} \joint\mathcal{C} =  \{\mathbf{BC}\dvtx\mathbf{B}
\in\mathcal{B}, \mathbf{C} \in\mathcal{C},  \mathbf{BC} \ne
\mathbf{0} \}.
\]
\end{definition}

Thus $\mathcal{B} \joint\mathcal{C} = \mathcal{C} \joint
\mathcal{B}$. Moreover, if $\mathcal{B}$ and $\mathcal{C}$ are also
both compatible with $\mathcal{D}$, then $\mathcal{B} \joint
\mathcal{C}$ is compatible with $\mathcal{D}$, $\mathcal{B}$ is
compatible with $\mathcal{C} \joint\mathcal{D}$, and $(\mathcal{B}
\joint\mathcal{C}) \joint\mathcal{D} = \mathcal{B} \joint
(\mathcal{C} \joint\mathcal{D})$. Hence if $\mathcal{B}_1, \ldots,
\mathcal{B}_m$ are pairwise compatible then there is no need for
parentheses in defining $\mathcal{B}_1 \joint\mathcal{B}_2 \joint
\cdots\joint\mathcal{B}_m$.
This decomposition could be
referred to as ``$\mathcal{B}_1$ combined with $\mathcal{B}_2$ combined
with $\cdots$ combined with $\mathcal{B}_m$.''
\begin{lemma}
\label{th:commute}
If $\mathbf{PQPRP}$ is symmetric for all $\mathbf{P}$ in $\mathcal{P}$,
all $\mathbf{Q}$ in $\mathcal{Q}$ and all $\mathbf{R}$ in~$\mathcal{R}$,
then $\mathcal{P} \vartriangleright\mathcal{Q}$
is compatible with $\mathcal{P} \vartriangleright\mathcal{R}$.
\end{lemma}
\begin{pf}
If
$\mathbf{PQPRP}$
is symmetric then
$\mathbf{PQPRP} = \mathbf{PRPQP}$.
Hence if
$\lambda_{\mathbf{PQ}} \times \lambda_{\mathbf{PR}} \ne0$ then
${(\mathbf{P} \vartriangleright\mathbf{Q})}{ (\mathbf{P}
\vartriangleright\mathbf{R})}
= \lambda_{\mathbf{PQ}}^{-1} \lambda_{\mathbf{PR}}^{-1}
\mathbf{PQPPRP} =\break \lambda_{\mathbf{PQ}}^{-1} \lambda_{\mathbf{PR}}^{-1}
\mathbf{PRPPQP} = (\mathbf{P} \vartriangleright\mathbf{R}) (\mathbf{P}
\vartriangleright
\mathbf{Q})$.
Thus if $\mathbf{P} \vartriangleright\mathbf{Q}$ is defined then it
commutes with $\mathbf{P}$ and with every $\mathbf{P}
\vartriangleright
\mathbf{R}$,
so it commutes with $\mathbf{P} \vdash\mathcal{R}$. Similarly,
if $\mathbf{P} \vartriangleright\mathbf{R}$ is defined then it
commutes with
$\mathbf{P} \vdash\mathcal{Q}$. Now the same argument shows that
$\mathbf{P} \vdash\mathcal{Q}$ commutes with $\mathbf{P} \vdash
\mathcal{R}$.
If $\mathbf{P}_i$ and $\mathbf{P}_j$ are different elements
of $\mathcal{P}$ then $\mathbf{P}_i \vartriangleright\mathbf{Q}$ and
$\mathbf{P}_i \vdash\mathcal{Q}$ commute with
$\mathbf{P}_j \vartriangleright\mathbf{R}$ and
$\mathbf{P}_j \vdash\mathcal{R}$, because all products are zero.
Hence $\mathcal{P} \vartriangleright\mathcal{Q}$ is compatible with
$\mathcal{P} \vartriangleright\mathcal{R}$.
\end{pf}
\begin{theorem}
\label{th:compat}
If the conditions in Lemma \ref{th:qrzero} are,
or condition (\ref{eq:coincond}) is,
satisfied, then $\mathcal{P} \vartriangleright\mathcal{Q}$
is compatible with $\mathcal{P} \vartriangleright\mathcal{R}$.
\end{theorem}
\begin{pf}
The first conditions imply that
$\mathbf{QPR} = \mathbf{0}$
or $\mathbf{P} = \mathbf{Q} = \mathbf{R} = \mathbf{P}_0$.
The second implies that $\mathbf{QPR} = \mathbf{0}$
or $\mathbf{PQP} = \lambda_{\mathbf{PQ}}\mathbf{P}$
or $\mathbf{PRP} = \lambda_{\mathbf{PR}}\mathbf{P}$.
In each case, $\mathbf{PQPRP}$
is symmetric, so Lemma \ref{th:commute} completes the proof.
\end{pf}

Thus the decomposition $(\mathcal{P} \vartriangleright\mathcal{Q})
\joint
(\mathcal{P} \vartriangleright\mathcal{R})$, which is symmetric in
$\mathcal{Q}$ and $\mathcal{R}$, can be used for coincident or
independent randomizations, or for
unrandomized-inclusive randomizations which satisfy
the conditions in Lemma \ref{th:qrzero}.
It is the same as decomposition
(\ref{eq:cdec}) for coincident randomizations
when condition (\ref{eq:coin}) holds, the same as
the decomposition in
Corollary \ref{eq:u2dec}
for unrandomized-inclusive randomizations when
the conditions in Lemma \ref{th:qrzero} hold,
and the same as decomposition
(\ref{eq:decindep}) for independent randomizations.
Condition (\ref{eq:coincond}) shows that, for a pair of coincident
randomizations, each idempotent in $(\mathcal{P} \vartriangleright
\mathcal{Q})
\joint(\mathcal{P} \vartriangleright\mathcal{R})$ has one of the following
forms: $\mathbf{P}$, $\mathbf{P} \vartriangleright\mathbf{Q}$,
$\mathbf{P}
\vartriangleright
\mathbf{R}$, $\mathbf{P} \vdash\mathcal{Q}$ or $\mathbf{P} \vdash
\mathcal{R}$.

If a pair of coincident randomizations does not satisfy
condition (\ref{eq:coin}), then it may be possible to refine
$\mathcal{R}$ to, say, $\mathcal{R}_2$ in such a way that
$\mathcal{R}_2$ is structure balanced
in relation to $\mathcal{P}
\vartriangleright\mathcal{Q}$, so that
the decomposition in Theorem \ref{th:RbalancedP}(b)
can be used.
It is possible if $\mathbf{R} = \mathbf{P}$ whenever
$\mathbf{P} \vartriangleright\mathbf{R} = \mathbf{P}$.
\setcounter{egg}{2}
\begin{egg}[(Continued)]
As already noted,
this example satisfies condition (\ref{eq:coin}), so
$\mathcal{P} \vartriangleright\mathcal{Q}$ is compatible with
$\mathcal{P} \vartriangleright\mathcal{R}$. Here
\begin{eqnarray*}
\mathcal{P} \vartriangleright\mathcal{Q} & = & \bigl\{\mathbf{P}_0
\vartriangleright\mathbf{Q}_0,  \mathbf{P}_{\mathrm{B}}
\vartriangleright\mathbf{Q}_{\mathrm{S}_1},  \mathbf{P}_{\mathrm
{P}[\mathrm{B}]} \vartriangleright\mathbf{Q}_{\mathrm{V}},
\mathbf{P}_{\mathrm{P}[\mathrm{B}]} \vartriangleright\mathbf
{Q}_{\mathrm{S}[\mathrm{V}] \vdash\mathrm{S}_1}  \bigr\}\\
& = &  \bigl\{\mathbf{P}_0,  \mathbf{P}_{\mathrm{B}},  \mathbf
{Q}_{\mathrm{V}},  \mathbf{Q}_{\mathrm{S}[\mathrm{V}] \vdash
\mathrm{S}_1} \bigr\}
\end{eqnarray*}
and
\begin{eqnarray*}
\mathcal{P} \vartriangleright\mathcal{R} &=&  \bigl\{\mathbf{P}_0
\vartriangleright\mathbf{R}_0,  \mathbf{P}_{\mathrm{B}}
\vartriangleright\mathbf{R}_{\mathrm{R}},  \mathbf{P}_{\mathrm
{B}} \vdash\mathcal{R},  \mathbf{P}_{\mathrm{P}[\mathrm{B}]}
\bigr\}\\
& = &
\bigl\{\mathbf{P}_0,  \mathbf{R}_{\mathrm{R}},  \mathbf
{P}_{\mathrm{B}} - \mathbf{R}_{\mathrm{R}},  \mathbf{P}_{\mathrm
{P}[\mathrm{B}]} \bigr\}.
\end{eqnarray*}

Then
\begin{eqnarray*}
(\mathcal{P} \vartriangleright\mathcal{Q}) \joint(\mathcal{P}
\vartriangleright
\mathcal{R})
& =&  \bigl\{\mathbf{P}_0^2,  \mathbf{P}_{\mathrm{B}}\mathbf
{R}_{\mathrm{R}},  \mathbf{P}_{\mathrm{B}}(\mathbf{P}_{\mathrm
{B}} - \mathbf{R}_{\mathrm{R}}),  \mathbf{Q}_{\mathrm{V}}\mathbf
{P}_{\mathrm{P}[\mathrm{B}]},  \mathbf{Q}_{\mathrm{S}[\mathrm
{V}] \vdash\mathrm{S}_1}\mathbf{P}_{\mathrm{P}[\mathrm{B}]}
\bigr\}
\nonumber\\
&=&  \bigl\{\mathbf{P}_0,  \mathbf{R}_{\mathrm{R}},  \mathbf
{P}_{\mathrm{B}} - \mathbf{R}_{\mathrm{R}},  \mathbf{Q}_{\mathrm
{V}},  \mathbf{Q}_{\mathrm{S}[\mathrm{V}] \vdash\mathrm
{S}_1} \bigr\}
\\
& = & (\mathcal{P} \vartriangleright\mathcal{Q}) \vartriangleright
\mathcal{R}.
\end{eqnarray*}
\end{egg}

\section{Double randomizations}
\label{s:double}

\begin{figure}

\includegraphics{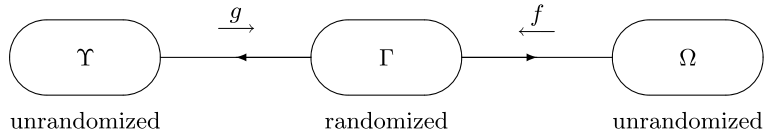}

\caption{Diagram of an experiment with double randomization.}
\label{f:Double}
\end{figure}

Double randomization is the one known type of two-from-one
randomizations. In an experiment with double randomization,
one set
of objects is randomized to two others; thus we could have
$\Gamma$ randomized to $\Upsilon$ and
to $\Omega$. We follow the
convention that the set of observational units is designated
as $\Omega$. Two functions are needed to encapsulate the results of
these randomizations, say ${f\colon\Omega\to\Gamma}$ and
${g\colon\Upsilon\to\Gamma}$. These two functions are
randomized independently using two different groups of
permutations. The set-up is shown in Figure \ref{f:Double}.

%

Now we obtain a subspace $V_{\Gamma}^f$ of $V_\Omega$
and a subspace $V_{\Gamma}^g$ of $V_{\Upsilon}$, both
isomorphic to $V_{\Gamma}$.
If $ |\Upsilon | =  |\Gamma |$ then
$V_{\Upsilon} =
V_{\Gamma}^g$, so we may effectively identify $V_{\Upsilon}$,
$V_{\Gamma}$ and $V_{\Gamma}^f$. If $ |\Upsilon | >
|\Gamma |$
then we cannot identify $V_{\Upsilon}$ with a subspace of $V_\Omega$
without further information explicitly assigning an element
of $\Upsilon$
to each observational unit in $\Omega$. This may not be possible (see,
e.g., Figure 28 in \cite{Brien06}). Thus we shall assume that
$ |\Upsilon | =  |\Gamma |$.

Associated with $\Omega$, $\Upsilon$ and $\Gamma$ are the decompositions
$\mathcal{P}$, $\mathcal{Q}$ and $\mathcal{R}$.
If $\mathcal{R}$ is structure balanced in relation to $\mathcal{Q}$
and $ |\Upsilon | =  |\Gamma |$, then
Lemma 4.2
in \cite{Brien09b} shows that $\mathcal{Q}\vartriangleright\mathcal
{R} =
\mathcal{R}$. Therefore it suffices to have $\mathcal{R}$ structure balanced
in relation to $\mathcal{P}$. Then the overall decomposition is
$\mathcal{P}
\vartriangleright\mathcal{R} = \mathcal{P} \vartriangleright
(\mathcal{Q} \vartriangleright
\mathcal{R})$, which must be done from right to left.
\begin{egg}[(An improperly replicated rotational grazing experiment)]
\label{eg:RotateGraze}
Example 8 in \cite{Brien06} is the rotational grazing trial shown in
Figure \ref{fig:improper}, with Cows substituted for Animals.
The double randomization of Availability
results in the assignment of Cows to Paddocks, the Cows assigned
to an Availability forming a single herd that is used to graze all
Paddocks with the same level of Availability. The sets of objects are
observational units, paddocks and treatments,
and the numbers of paddocks and treatments are equal, as required.
The Hasse diagrams for
treatments and observational units are like
the middle diagram in Figure \ref{f:HassePlant};
that for paddocks is trivial.

%
\begin{figure}[b]

\includegraphics{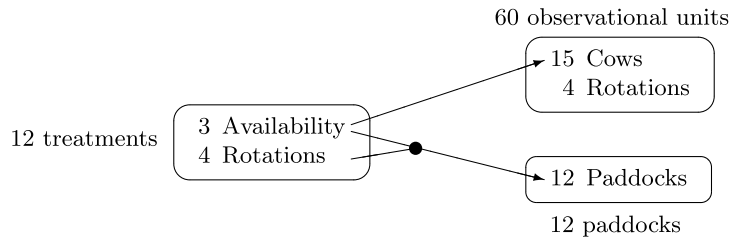}

\caption{Double randomizations in Example \protect\ref
{eg:RotateGraze}: treatments
are randomized to both observational units and paddocks.}
\label{fig:improper}
\end{figure}

The structures on observational units,
paddocks and treatments are
$\mathcal{P} = \{\mathbf{P}_0, \mathbf{P}_{\mathrm
{C}}, \mathbf{P}_{\mathrm{R}}, \mathbf{P}_{\mathrm{C}\#\mathrm
{R}} \}$,
$\mathcal{Q} =  \{\mathbf{Q}_0, \mathbf{Q}_{\mathrm{P}}
\}$ and
$\mathcal{R} =  \{\mathbf{R}_0, \mathbf{R}_{\mathrm{A}},
\mathbf{R}_{\mathrm{R}}, \mathbf{R}_{\mathrm{A}\#\mathrm
{R}} \}$, respectively.
This leads to the decomposition
$\mathcal{P} \vartriangleright(\mathcal{Q} \vartriangleright
\mathcal{R})$ in
Table \ref{tab:ANOVARotateGraze}. It shows that
there are no residual degrees of freedom for testing
any treatment differences---hence the experiment being dubbed \textit{improperly replicated}.

In this case,
\[
\mathcal{Q} \vartriangleright\mathcal{R} = \mathcal{R} =
 \{\mathbf{Q}_0\vartriangleright\mathbf{R}_0, \mathbf
{Q}_{\mathrm{P}} \vartriangleright\mathbf{R}_{\mathrm{A}},  \mathbf
{Q}_{\mathrm{P}} \vartriangleright\mathbf{R}_{\mathrm{R}},  \mathbf
{Q}_{\mathrm{P}} \vartriangleright\mathbf{R}_{\mathrm{A} \# \mathrm
{R}} \}
\]
with
\begin{eqnarray*}
\mathbf{Q}_0\vartriangleright\mathbf{R}_0 &=& \mathbf{R}_0,\qquad
\mathbf{Q}_{\mathrm{P}} \vartriangleright\mathbf{R}_{\mathrm{A}} =
\mathbf
{R}_{\mathrm{A}},\\
\mathbf{Q}_{\mathrm{P}} \vartriangleright\mathbf{R}_{\mathrm{R}} &=&
\mathbf
{R}_{\mathrm{R}},\qquad
\mathbf{Q}_{\mathrm{P}} \vartriangleright\mathbf{R}_{\mathrm{A} \#
\mathrm{R}} =
\mathbf{R}_{\mathrm{A} \#\mathrm{R}}.
\end{eqnarray*}

Also, $\mathbf{PR}$ is equal to either $\mathbf{R}$ or $\mathbf{0}$
for all $\mathbf{P} \in\mathcal{P}$ and all $\mathbf{R} \in
\mathcal{R}$.
That is, $\mathcal{R}$ is orthogonal in relation to $\mathcal{P}$.
Therefore the complete decomposition for the experiment is
\[
\mathcal{P} \vartriangleright(\mathcal{Q} \vartriangleright\mathcal
{R}) =  \left\{
\begin{array}{l}
\mathbf{P}_0\vartriangleright(\mathbf
{Q}_0\vartriangleright\mathbf{R}_0),  \mathbf{P}_{\mathrm{C}}
\vartriangleright(\mathbf{Q}_{\mathrm{P}} \vartriangleright\mathbf
{R}_{\mathrm{A}}),  \mathbf{P}_{\mathrm{C}} \vdash(\mathcal{Q}
\vartriangleright\mathcal{R}) , \\
\mathbf{P}_{\mathrm{R}}
\vartriangleright(\mathbf{Q}_{\mathrm{P}} \vartriangleright\mathbf
{R}_{\mathrm{R}}),  \mathbf{P}_{\mathrm{C} \#\mathrm{R}}
\vartriangleright(\mathbf{Q}_{\mathrm{P}} \vartriangleright\mathbf
{R}_{\mathrm{A} \# \mathrm{R}}),  \mathbf{P}_{\mathrm{C} \#\mathrm
{R}} \vdash(\mathcal{Q} \vartriangleright\mathcal{R})
\end{array}
\right\}
\]
with
\begin{eqnarray*}
\mathbf{P}_0\vartriangleright\mathbf{Q}_0\vartriangleright\mathbf
{R}_0 &=&
\mathbf{P}_0,\qquad
\mathbf{P}_{\mathrm{C}} \vartriangleright(\mathbf{Q}_{\mathrm{P}}
\vartriangleright
\mathbf{R}_{\mathrm{A}}) =
\mathbf{Q}_{\mathrm{P}} \vartriangleright\mathbf{R}_{\mathrm{A}} =
\mathbf
{R}_{\mathrm{A}}, \\
\mathbf{P}_{\mathrm{C}} \vdash(\mathcal{Q} \vartriangleright
\mathcal{R}) &=&
\mathbf{P}_{\mathrm{C}} -(\mathbf{Q}_{\mathrm{P}} \vartriangleright
\mathbf
{R}_{\mathrm{A}})
= \mathbf{P}_{\mathrm{C}}- \mathbf{R}_{\mathrm{A}}, \\
\mathbf{P}_{\mathrm{R}} \vartriangleright(\mathbf{Q}_{\mathrm{P}}
\vartriangleright
\mathbf{R}_{\mathrm{R}}) &=&
\mathbf{Q}_{\mathrm{P}} \vartriangleright\mathbf{R}_{\mathrm{R}} =
\mathbf
{R}_{\mathrm{R}}, \\
\mathbf{P}_{\mathrm{C} \#\mathrm{R}} \vartriangleright
(\mathbf{Q}_{\mathrm{P}} \vartriangleright\mathbf{R}_{\mathrm{A} \#
\mathrm{R}}) &=&
\mathbf{Q}_{\mathrm{P}} \vartriangleright\mathbf{R}_{\mathrm{A} \#
\mathrm{R}} =
\mathbf{R}_{\mathrm{A} \#\mathrm{R}}, \\
\mathbf{P}_{\mathrm{C} \#\mathrm{R}} \vdash(\mathcal{Q}
\vartriangleright\mathcal{R}) &=&
\mathbf{P}_{\mathrm{C} \#\mathrm{R}} -
(\mathbf{Q}_{\mathrm{P}} \vartriangleright\mathbf{R}_{\mathrm{A} \#
\mathrm{R}})
= \mathbf{P}_{\mathrm{C} \#\mathrm{R}} -
\mathbf{R}_{\mathrm{A} \#\mathrm{R}}.
\end{eqnarray*}

%
\begin{table}
\caption{Decomposition table
for Example \protect\ref{eg:RotateGraze}}\label{tab:ANOVARotateGraze}
\begin{tabular*}{\tablewidth}{@{\extracolsep{4in minus 4in}}lrclrclr@{}}
\hline
\multicolumn{2}{@{}c}{\textbf{observational units tier}} &&
\multicolumn{2}{c}{\textbf{paddocks tier}} &&
\multicolumn{2}{c}{\textbf{treatments tier}}\\[-4pt]
\multicolumn{2}{@{}l}{\hrulefill} && \multicolumn
{2}{c}{\hrulefill} && \multicolumn{2}{c@{}}{\hrulefill} \\
\textbf{source} & \textbf{d.f.} && \textbf{source} & \textbf{d.f.}
&& \textbf{source} & \textbf{d.f.} \\
\hline
Mean & \phantom{0}1 && Mean & \phantom{0}1 && Mean & 1\\
\hline
Cows & 14 && Paddocks & \phantom{0}2 && Availability & 2\\[-4pt]
&&& \multicolumn{5}{c@{}}{\hrulefill}\\
& & & Residual & 12\\
\hline
Rotations & \phantom{0}3 && Paddocks & \phantom{0}3 && Rotations & 3\\
\hline
$\mathrm{Cows}\,\#\,\mathrm{Rotations}$ & 42 &&
Paddocks & \phantom{0}6 && $\mathrm{Availability} \,\#\,\mathrm{Rotations}$
&6\\
[-4pt]
&&& \multicolumn{5}{c@{}}{\hrulefill}\\
& & & Residual & 36 \\
\hline
\end{tabular*}
\end{table}

In \cite{Brien06}
this example was redone as a case of
randomized-inclusive
randomization, using two pseudofactors $\mathrm{P}_{\mathrm{A}}$ and
$\mathrm{P}_{\mathrm{R}}$ for Paddocks, aliased with Availability and
Rotations, respectively. These are required if $\mathcal{Q}$ itself is
to be structure balanced
in relation to $\mathcal{P}$, giving a
decomposition from left to right like the one in Section 6
in \cite{Brien09b}.
\end{egg}

\section{Summary}
\label{s:summary2}

We have shown in \cite{Brien09b} and here that,
under structure balance, the
six different types of multiple randomization
identified in \cite{Brien06} all lead to orthogonal
decompositions of $V_\Omega$ using some of the following
idempotents: $\mathbf{P}$, $\mathbf{P} \vartriangleright\mathbf{Q}$,
$\mathbf{P} \vartriangleright\mathbf{R}$, $(\mathbf{P}
\vartriangleright
\mathbf{Q}) \vartriangleright\mathbf{R}$, $\mathbf{P}
\vartriangleright(\mathbf{Q}
\vartriangleright\mathbf{R})$, $\mathbf{P} \vdash\mathcal{Q}$,
$\mathbf{P} \vdash\mathcal{R}$, $(\mathbf{P} \vartriangleright
\mathbf{Q})
\vdash\mathcal{R}$, $\mathbf{P} \vartriangleright(\mathbf{Q}
\vdash
\mathcal{R})$, $(\mathbf{P} \vdash\mathcal{Q}) \vartriangleright
\mathbf{R}$ and $(\mathbf{P} \vdash\mathcal{Q}) \vdash
\mathcal{R}$. The differences between the different multiple
randomizations lead to differences in the
reduced forms for these elements and in the efficiency factors.

\subsection*{Composed randomizations}
If each design is structure balanced
then so is the composite; the decompositions
$\mathcal{P} \vartriangleright(\mathcal{Q} \vartriangleright
\mathcal{R})$ and
$(\mathcal{P} \vartriangleright\mathcal{Q}) \vartriangleright
\mathcal{R}$ are
equal, and so the decomposition may be done in either order; and there
are no
idempotents of the form $(\mathbf{P} \vdash\mathcal{Q})
\vartriangleright
\mathbf{R}$
or $(\mathbf{P} \vdash\mathcal{Q}) \vdash\mathcal{R}$.

\subsection*{Randomized-inclusive randomizations}
The structures $\mathcal{Q}_1$ and $\mathcal{R}_1$ for design~1
are refined to $\mathcal{Q}$ and $\mathcal{R}$ using the
pseudofactors that are necessary for the second randomization,
and then the results
are the same as for composed randomizations.

\subsection*{Unrandomized-inclusive randomizations}
We must have $\mathcal{R}$ structure balanced
in relation to
$\mathcal{P} \vartriangleright\mathcal{Q}$; use the decomposition
$(\mathcal{P} \vartriangleright\mathcal{Q}) \vartriangleright
\mathcal{R}$, which
is done from left to right; if
the conditions in Lemma \ref{th:qrzero} hold
then there are no
idempotents of the form $(\mathbf{P} \vartriangleright\mathbf{Q})
\vartriangleright
\mathbf{R}$
apart from the Mean,
nor any of the form $(\mathbf{P} \vartriangleright\mathbf{Q}) \vdash
\mathcal{R}$,
the decomposition
$\mathcal{P} \vartriangleright\mathcal{Q}$ is compatible with
$\mathcal{P} \vartriangleright\mathcal{R}$,
and $(\mathcal{P} \vartriangleright\mathcal{Q})
\vartriangleright\mathcal{R} = (\mathcal{P} \vartriangleright
\mathcal{Q}) \joint
(\mathcal{P}
\vartriangleright\mathcal{R})$.

\subsection*{Independent randomizations}
The conditions in Lemma \ref{th:qrzero}
must hold; if both designs are structure balanced
then each remains structure balanced
after the other has been taken into
account; $\mathcal{P} \vartriangleright\mathcal{Q}$ is compatible with
$\mathcal{P}
\vartriangleright\mathcal{R}$; the decompositions $(\mathcal{P}
\vartriangleright
\mathcal{Q})
\vartriangleright\mathcal{R}$, $(\mathcal{P} \vartriangleright
\mathcal{R}) \vartriangleright
\mathcal{Q}$ and $(\mathcal{P} \vartriangleright\mathcal{Q}) \joint
(\mathcal{P}
\vartriangleright\mathcal{R})$ are equal; and there are no
idempotents of the form $(\mathbf{P} \vartriangleright\mathbf{Q})
\vartriangleright
\mathbf{R}$
apart from the Mean, nor any of the form
$(\mathbf{P} \vartriangleright\mathbf{Q}) \vdash\mathcal{R}$.

\subsection*{Coincident randomizations}
Condition (\ref{eq:coincond}) must hold;
$\mathcal{P} \vartriangleright\mathcal{Q}$ is compatible with
$\mathcal{P}
\vartriangleright
\mathcal{R}$; use the decomposition $(\mathcal{P} \vartriangleright
\mathcal
{Q}) \joint
(\mathcal{P} \vartriangleright\mathcal{R})$,
whose idempotents have the form $\mathbf{P}$, $\mathbf{P}
\vartriangleright
\mathbf{Q}$, $\mathbf{P} \vartriangleright\mathbf{R}$, $\mathbf{P}
\vdash
\mathcal{Q}$ or $\mathbf{P} \vdash\mathcal{R}$;
if condition (\ref{eq:coin}) holds, this is the same as the
decomposition $(\mathcal{P} \vartriangleright\mathcal{Q})
\vartriangleright
\mathcal{R}$, which is done from left to right; otherwise, there
may be a refinement of $\mathcal{R}$ giving a left-to-right
decomposition.

\subsection*{Double randomizations}

We require that $ |\Upsilon | =  |\Gamma |$ and that
$\mathcal{R}$ be
structure balanced in relation to both $\mathcal{Q}$ and $\mathcal{P}$,
so that the decomposition is $ \mathcal{P} \vartriangleright\mathcal
{R} =
\mathcal{P} \vartriangleright(\mathcal{Q} \vartriangleright\mathcal
{R})$, which is
done from
right to left.
It appears that they can also be formulated as
randomized-inclusive randomizations using pseudofactors to
refine $\mathcal{Q}$ to $\mathcal{Q}_2$ for which the
left-to-right decomposition $(\mathcal{P} \vartriangleright\mathcal{Q}_2)
\vartriangleright\mathcal{R}$ is correct.

\section{Structure-balanced experiments with four or more tiers}
\label{s:fourII}

Each experiment in Sections \ref{s:uincl}--\ref{s:double} involves
only one type of multiple randomization, and so involves three
tiers and three structures. However, multitiered experiments are
not limited to this configuration. Examples 12--14 in \cite{Brien06}
each have four tiers and involve more than one
type of multiple randomization. In general, there is the set of
observational units, $\Omega$, and each randomization adds another
set of objects with its associated tier.

Section 7 in \cite{Brien09b} shows how to deal with
three or more randomizations which follow each other in a chain.
Mixtures of other types of multiple
randomization should be amenable to successive decompositions of the
sort summarized in Section \ref{s:summary2},
so long as they are handled in the correct order.
Thus we can use a recursive
procedure in which each new structure refines the decomposition
of $V_\Omega$ obtained using structures accounted for previously.
All that is required is that each successive structure
should be structure balanced
in relation to the previous decomposition.

%
\begin{figure}[b]

\includegraphics{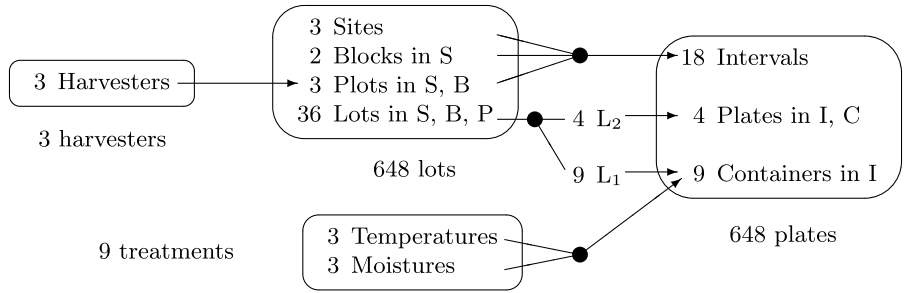}

\caption{Composed and coincident randomizations in
Example \protect\ref{eg:corn}: harvesters are randomized to plots; lots
of grain are sampled from each plot and then randomized to
plates; and treatments are randomized to plates; $\mathrm{S}$,
$\mathrm{B}$, $\mathrm{P}$, $\mathrm{I}$, $\mathrm{C}$ denote Sites,
Blocks, Plots, Intervals, Containers, respectively.}
\label{fig:corn}
\end{figure}

One class of experiments with both two--one randomizations and chain
randomizations consists of multiphase experiments in which different
treatment factors are applied in different phases, as the following
example demonstrates.
\begin{egg}[(A two-phase corn seed germination experiment)]
\label{eg:corn}
Example 12 of \cite{Brien06}
has the four tiers shown in Figure \ref{fig:corn}.
Here we have taken the opportunity to correct the diagram given
in \cite{Brien06}.
The $36$ Lots of grain within each Plot should be completely randomized to
$\mathrm{Plates}\wedge\mathrm{Containers}$ within each Interval.
This will not be achieved by permuting Containers within Intervals and
Plates within $\mathrm{Intervals}\wedge\mathrm{Containers}$, as implied
in the rightmost panel of Figure \ref{fig:corn}.
We introduce pseudofactors $\mathrm{L}_1$ and $\mathrm{L}_2$ for
Lots, with nine and four levels, respectively,
like the pseudofactors for Seedlings in Example \ref{eg:PlantExp}. The
$36$ Lots must be randomly allocated to the combinations of levels of
$\mathrm{L}_1$ and $\mathrm{L}_2$, independently within each level of
$\mathrm{Sites} \wedge\mathrm{Blocks} \wedge\mathrm{Plots}$, so
that neither
pseudofactor corresponds to any inherent source of variation.

At each randomization, an orthogonal design is used, so there is
no difficulty in constructing the decomposition in
Table \ref{tab:corn}.
Here $\mathrm{L}_1 [\mathrm{P} \wedge\mathrm{B} \wedge\mathrm{S}]$
is the part of the source $\mathrm{Lots} [\mathrm{P} \wedge\mathrm
{B} \wedge\mathrm{S}]$ which is confounded with $\mathrm
{Containers} [\mathrm{I}]$.
The source $\mathrm{Lots} [\mathrm{P} \wedge\mathrm{B} \wedge
\mathrm{S}]_{\vdash}$
is the part of $\mathrm{Lots} [\mathrm{P} \wedge\mathrm{B} \wedge
\mathrm{S}]$ which is orthogonal to
$\mathrm{L}_1 [\mathrm{P} \wedge\mathrm{B} \wedge\mathrm{S}]$:
it is confounded with $\mathrm{Plates} [\mathrm{C}\wedge\mathrm{I}]$.

%
\begin{table}
\tabcolsep=0pt
\caption{Decomposition table for Example
\protect\ref{eg:corn}}\label{tab:corn}
\begin{tabular*}{\tablewidth}{@{\extracolsep{4in minus 4in}}lrlrlrlr@{}}
\hline
\multicolumn{2}{@{}c}{\textbf{plates tier}} &
\multicolumn{2}{c}{\textbf{lots tier}} &
\multicolumn{2}{c}{\textbf{harvesters tier}} &
\multicolumn{2}{c@{}}{\textbf{treatments tier}}\\[-4pt]
\multicolumn{2}{@{}l}{\hrulefill} & \multicolumn{2}{c}{\hrulefill} & \multicolumn{2}{c}{\hrulefill} &
\multicolumn{2}{c@{}}{\hrulefill} \\
\textbf{source} & \multicolumn{1}{c}{\textbf{d.f.}}
& \multicolumn{1}{l}{\textbf{source}} &
\multicolumn{1}{l}{\textbf{d.f.}}
& \multicolumn{1}{l}{\textbf{source}} & \multicolumn{1}{l}{\textbf{d.f.}}
& \multicolumn{1}{l}{\textbf{source}} & \multicolumn{1}{l@{}}{\textbf{d.f.}} \\
\hline
Mean & \phantom{00}1 & Mean & \phantom{00}1 & Mean & \phantom{0}1 & Mean & \phantom{00}1\\
\hline
Intervals& \phantom{0}17 & Sites & \phantom{00}2\\[-4pt]
&& \multicolumn{6}{c@{}}{\hrulefill}\\
& & $\mathrm{Blocks}[\mathrm{S}]$ & \phantom{00}3\\[-4pt]
&& \multicolumn{6}{c@{}}{\hrulefill}\\
& & $\mathrm{Plots}[\mathrm{B}\wedge\mathrm{S}]$ &
\phantom{0}12 & Harvesters &\phantom{0}2\\
& & & & Residual & 10\\
\hline
$\mathrm{Containers} [\mathrm{I}]$ & 144 &
$\mathrm{L}_1 [\mathrm{P} \wedge\mathrm{B} \wedge\mathrm{S}]$
& 144 & & &Temperature & \phantom{00}2\\
& & & & & & Moistures& \phantom{00}2\\
& & & & & & $\mathrm{T} \,\#\,\mathrm{M}$ & \phantom{00}4\\
& & & & & & Residual & 136\\
\hline
$\mathrm{Plates}[\mathrm{C}\wedge\mathrm{I}]$
& 486 &
$\mathrm{Lots} [\mathrm{P} \wedge\mathrm{B} \wedge\mathrm
{S}]_{\vdash}$ & 486\\
\hline
\end{tabular*}
\end{table}

Bailey \cite{baileytb} suggests an analysis for this example which we
reproduce in the first three columns of
Table \ref{tab:Baileycorn}(a). In this, the $3$-level factors
Temperature and Moisture have been combined into a single $9$-level
Treatment factor, the intertier interactions \cite{Brien06} of
Sites, Harvesters and Treatments have been included, and the notation
$\times$ is used in place of $\#$. We cannot be sure, but it is
plausible that he
based this decomposition on the crossing and nesting relationships
summarized in the formula
%
%
\begin{equation}\label{eq:tbb}
\bigl(\mathrm{T} *\mathrm{H} *(\mathrm{S}/\mathrm{B})\bigr)
/\mathrm{Q},
\end{equation}
where $\mathrm{T}$, $\mathrm{H}$, $\mathrm{S}$, $\mathrm{B}$ and
$\mathrm{Q}$ represent factors for Treatments, Harvesters, Sites,
Blocks and Plates, with $9$, $3$, $3$, $2$ and $4$ levels, respectively.
The sources derived from this are in the final column of
Table \ref{tab:Baileycorn}(a), with degrees of freedom matching those in
the preceding column.

%
\begin{table}
\caption{Skeleton analysis-of-variance tables
for Example \protect\ref{eg:corn}\textup{(a)} given by
Bailey \protect\cite{baileytb} and
\textup{(b)} from Table \protect\ref{tab:corn}
with intertier interactions added}
\label{tab:Baileycorn}
\begin{tabular*}{\tablewidth}{@{\extracolsep{\fill}}llrl@{}}
(a) \\
\hline
& \textbf{source} & \textbf{d.f.} & \textbf{source from (\ref{eq:tbb})} \\
\hline
\textbf{Phase I: field study} & Site & 2 & $\mathrm{S}$\\
& Experimental error (a) & 3  & $\mathrm{B}[\mathrm{S}]$\\
& Harvester & 2 & $\mathrm{H}$\\
&$\mathrm{Harvester} \times\mathrm{Site}$ & 4 & $\mathrm{H}\,\#\,
\mathrm{S}$\\
& Experimental error (b) & 6 &  $\mathrm{H}\,\#\,\mathrm{B}[\mathrm
{S}]$\\
[4pt]
\textbf{Phase II: laboratory study} & Treatment & 8 &  $\mathrm{T}$\\
& $\mathrm{Treatment} \times\mathrm{Site}$ & 16 &
$\mathrm{T} \,\#\,\mathrm{S}$\\
& Experimental error (c) & 24 &  $\mathrm{T} \,\#\,\mathrm
{B}[\mathrm{S}]$\\
&
$\mathrm{Treatment} \times\mathrm{Harvester}$ & 16 &
$\mathrm{T} \,\#\,\mathrm{H}$\\
& $\mathrm{Treatment} \times\mathrm{Harvester} \times\mathrm{Site}$ & 32
& $\mathrm{T} \,\#\,\mathrm{H} \,\#\,\mathrm{S}$\\
&Experimental error (d) & 48 &
$\mathrm{T} \,\#\,\mathrm{H} \,\#\,\mathrm{B}[\mathrm{S}]$\\
& Residual & 486&  $\mathrm{Q}[\mathrm{T} \wedge\mathrm{H} \wedge
\mathrm{S} \wedge\mathrm{B}]$\\
\hline
\end{tabular*}
\begin{tabular*}{\tablewidth}{@{\extracolsep{4in minus 4in}}lrlrlrlr@{}}
(b) \\
\hline
\multicolumn{2}{@{}c}{\textbf{plates tier}} &
\multicolumn{2}{c}{\textbf{lots tier}} &
\multicolumn{2}{c}{\textbf{harvesters tier}} &
\multicolumn{2}{c@{}}{\textbf{treatments tier}}\\[-4pt]
\multicolumn{2}{@{}l}{\hrulefill} & \multicolumn{2}{c}{\hrulefill} & \multicolumn{2}{c}{\hrulefill} &
\multicolumn{2}{c@{}}{\hrulefill} \\
\textbf{source} & \textbf{d.f.} &  \textbf{source} & \textbf{d.f.}
& \textbf{source} & \textbf{d.f.} & \textbf{source} &
\textbf{d.f.} \\
\hline
Mean & \phantom{00}1 & Mean & \phantom{00}1 & Mean & 1 & Mean & \phantom{0}1 \\
\hline
Intervals& \phantom{0}17 & Sites & \phantom{00}2\\
\hline
& & $\mathrm{Blocks}[\mathrm{S}]$ & \phantom{00}3  \\[-4pt]
&& \multicolumn{6}{c@{}}{\hrulefill}\\
& & $\mathrm{Plots}[\mathrm{B}\wedge\mathrm{S}]$ &
\phantom{0}12 & Harvesters & 2 \\
& & & & $\mathrm{H} \,\#\,\mathrm{S}$ & 4 \\
& & & & Residual & 6 \\
\hline
$\mathrm{Containers} [\mathrm{I}]$ & 144 &
$\mathrm{Lots} [\mathrm{P} \wedge\mathrm{B} \wedge\mathrm{S}]_1$
& 144 & & & Treatments & \phantom{0}8 \\
& & & & & & $\mathrm{T} \,\#\,\mathrm{S}$ & 16 \\
& & & & & & $\mathrm{T} \,\#\,\mathrm{H}$ & 16 \\
& & & & & & $\mathrm{T} \,\#\,\mathrm{H} \,\#\,\mathrm{S}$ &
32 \\
& & & & & & Residual & 72\\
\hline
$\mathrm{Plates}[\mathrm{C}\wedge\mathrm{I}]$
& 486 &
$\mathrm{Lots} [\mathrm{P} \wedge\mathrm{B} \wedge\mathrm
{S}]_{\vdash}$ & 486\\
\hline
\end{tabular*}
\end{table}

Revision of Table \ref{tab:corn} along similar lines, and with
pseudosources replaced with actual sources, yields the skeleton
analysis-of-variance table in Table \ref{tab:Baileycorn}(b). Note
that, given Step 4 in Table 1 of \cite{Brien83}, an intertier
interaction will generally occur in the right-most tier that contains
a main effect in the interaction. Table \ref{tab:Baileycorn}(b)
differs from Table \ref{tab:Baileycorn}(a) in the following ways.

\begin{enumerate}
\item
The rationale for the sources in Table \ref{tab:Baileycorn}(a) is
unclear. We had to reverse-engineer it by producing formula (\ref{eq:tbb}).
On the other hand, the sources in Table \ref{tab:Baileycorn}(b) are
based on the relationships between factors within each tier and on the
confounding between sources from different tiers.
\item
Table \ref{tab:Baileycorn}(a) does not show,
as Table \ref{tab:Baileycorn}(b) does,
the successive decomposition of the vector space indexed by the
observational units. The impression given is that there is
a set of sources that arise from the field phase and another set that
arises from the laboratory phase.
\item
Table \ref{tab:Baileycorn}(a) has four sources called ``experimental
error'' and does not mention plates, containers,
intervals, blocks, plots or lots. Hence, there is no indication
of the sources of error variation. By contrast,
each source called ``Residual'' in Table \ref{tab:Baileycorn}(b) is
unambiguously identified; and the labelling
shows that all terms are affected by variation from both phases.
For example, the Residual for $\mathrm{Plots}[\mathrm{B}\wedge
\mathrm{S}]$, labelled Experimental error (b) in
Table \ref{tab:Baileycorn}(a), clearly arises from variability
associated with Plots within the Sites-Blocks combinations and
variability associated with Intervals. Similarly, it can be seen from
Table \ref{tab:Baileycorn}(b) that the Residual in
Table \ref{tab:Baileycorn}(a) arises
from variability associated with Plates and Lots.
\item
As discussed in \cite{Brien06}, Section 7.1, the usual default is that
there are no
intertier interactions because such inclusions would mean that the
analysis cannot be justified by the randomization used. It parallels
the assumption of
unit-treatment additivity in single-randomization experiments.
The approach
using Table \ref{tab:corn} forces the statistician to to consult the
researcher about
whether intertier interactions should be included, and, if so, to
justify them.
Tables \ref{tab:Baileycorn}(a) and (b) include the intertier interactions
of Sites, Harvesters and Treatments, which suggests that it is
anticipated that Harvesters
and Treatments will perform differently at different Sites.

\item
Even with the addition of intertier interactions, the decompositions in
Tables~\ref{tab:Baileycorn}(a) and (b) are not equivalent, and so
neither are the mixed models underlying them. Experimental errors (c)
and (d) from Table \ref{tab:Baileycorn}(a) are combined into the
Residual with $72$ degrees of freedom for $\mathrm{Lots} [\mathrm{P}
\wedge\mathrm{B} \wedge\mathrm{S}]_1$ in Table \ref
{tab:Baileycorn}(b). To justify an analysis based on Table \ref
{tab:Baileycorn}(a), one would need to argue that
unit-treatment interaction of Treatments with Blocks within Sites
can be anticipated in this experiment.
\end{enumerate}
\end{egg}

\section{Discussion}
\label{s:discuss}

\subsection{Implications of incoherent unrandomized-inclusive randomizations}

The phenomenon of incoherent unrandomized-inclusive randomizations is
described in~\cite{Brien06}, Section 5.2.1.
Essentially, when there has been a randomization to factors
that are crossed, one or more of these factors become nested in the
second randomization.

Consider the cherry rootstock experiment in Example \ref{eg:SuperStruct}.
The trees tier gives an orthogonal decomposition of $V_\Omega$ into
sources Mean, Blocks and $\mathrm{Trees} [\mathrm{Blocks}]$
of dimensions $1$, $2$ and $27$, respectively, in the left-hand
column of Table \ref{tab:ANOVASuper}. Similarly, the rootstocks tier
decomposes $V_\Upsilon$ into sources Mean and Rootstocks
of dimensions $1$ and $9$. The result of the first randomization is to
make the Mean sources equal and to place the Rootstocks source
inside $\mathrm{Trees} [\mathrm{Blocks}]$,
thus giving the finer decomposition
of $V_\Omega$ shown in the middle column of Table \ref{tab:ANOVASuper}.

The result of the
unrandomized-inclusive randomization should be to further
decompose the decomposition resulting from the first two tiers.
In the extended Youden square,
the source Viruses is
orthogonal to Blocks but
partially confounded with Rootstocks, so the Viruses source defines
the decomposition in the right-hand column in
Table \ref{tab:ANOVASuper}. That is, the source
Viruses further decomposes the sources
Rootstocks and the Residual for $\mathrm{Trees}
[\mathrm{Blocks}]$, as required.

In \cite{Brien06}
we discussed the possibility that the designer of the
superimposed experiment ignores the inherent crossing of the factors
Blocks and
Rootstocks and randomizes Viruses to Blocks in Rootstocks in a
balanced\break incomplete-block design. Then the randomizations are
incoherent. The permutation group for the second randomization does
not preserve the structure arising from the first two tiers, exhibited
by the two left-most columns in Table \ref{tab:ANOVASuper}. We can see
immediately that this randomization is senseless because it destroys
the Blocks subspace
preserved by the first randomization. This randomization
might have some appeal if no block effects had been detected during
the 20 years of the original experiment, but then the analysis of the
second experiment would be based on
an assumed model rather than on the intratier structures.

%
\begin{figure}[b]

\includegraphics{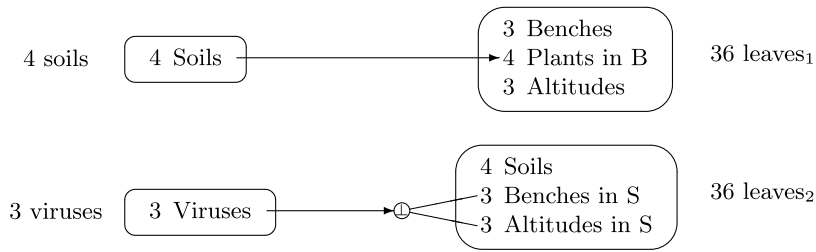}

\caption{Incoherent randomizations in Example \protect\ref{eg:cc}: both
soils and viruses are randomized to leaves, but with
different structures on leaves; $\mathrm{B}$ denotes Benches;
$\mathrm{S}$ denotes Soils.}
\label{fig:cc}
\end{figure}

%
\begin{table}[b]
\caption{Attempted decomposition table for
Example \protect\ref{eg:cc}}\label{tab:cc}
\begin{tabular*}{\tablewidth}{@{\extracolsep{4in minus 4in}}lrlrlrclr@{}}
\hline
\multicolumn{2}{@{}c}{\textbf{leaves$_1$ tier}} &
\multicolumn{2}{c}{\textbf{soils tier}} &
\multicolumn{2}{c}{\textbf{leaves$_2$ tier}\tabnoteref{t1}} &
\multicolumn{3}{c@{}}{\textbf{viruses tier}}\\[-4pt]
\multicolumn{2}{@{}l}{\hrulefill} & \multicolumn{2}{c}{\hrulefill} & \multicolumn{2}{c}{\hrulefill} &
\multicolumn{3}{c@{}}{\hrulefill} \\
\textbf{source} & \textbf{d.f.} & \textbf{source} &
\textbf{d.f.} & \textbf{source} & \textbf{d.f.} &
\textbf{eff.} & \textbf{source} & \textbf{d.f.} \\
\hline
Mean & 1 & Mean & 1 & Mean & 1 && Mean & 1\\
\hline
Benches & 2  & & & $\mathrm{B}[\mathrm{S}]_{\mathrm{B}}$ &
2  \\
\hline
$\mathrm{Plants}[\mathrm{B}]$ & 9 &  Soils & 3 &  Soils &3
\\[-4pt]
&& \multicolumn{7}{c@{}}{\hrulefill}\\
& &Residual &6 & $\mathrm{B}[\mathrm{S}]_{\vdash}$ & 6  \\
\hline
Altitudes & 2 & & &$\mathrm{A}[\mathrm{S}]_{\mathrm{A}}$
& 2 \\
\hline
A$\,\#\,$B & 4 & &&
$\mathrm{A}\,\#\,\mathrm{B}[\mathrm{S}]_{\mathrm{A}\,\#\,
\mathrm{B}}$
& 4 &? & Viruses\tabnoteref{t2} & 2 \\
& & & & & & ? & V$\,\#\,$S\tabnoteref{t2} & 2 \\
\hline
$\mathrm{A}\,\#\,\mathrm{P}[\mathrm{B}]$ & 18 && &
$\mathrm{A}[\mathrm{S}]_\vdash$ & 6  \\[-4pt]
&& \multicolumn{7}{c@{}}{\hrulefill}\\
& & & &$\mathrm{A}\,\#\,\mathrm{B}[\mathrm{S}]_\vdash
$ &
12  & ?& Viruses\tabnoteref{t2} & 2 \\
 & & & & & & ?& V$\,\#\,$S\tabnoteref{t2} & 6 \\
 & & & & & & & Residual & 4 \\
\hline
\end{tabular*}
\tabnotetext[\dag]{t1}{The subscripts on the sources from this tier indicate that they
are the part of the source associated with the subscripted source in
the first tier, and the subscript ``$\vdash$'' that this is the part of
the source orthogonal to all previous parts.}
\tabnotetext[\ddag]{t2}{The partial confounding of Viruses and V$\,\#\,$S may not have
first-order balance.}
\end{table}

Other examples of incoherent unrandomized-inclusive randomizations are more
complicated, and perhaps less easily detected. One is the design
proposed by several authors for a split-plot experiment in which the subplot
treatments are to be assigned using a row-column
design. Example \ref{eg:cc} illustrates how consideration of the
decomposition table for the proposed design facilitates the design
process and helps the detection of incoherence.
\begin{egg}[(Split-plots in a row--column design)]\label{eg:cc}
Example 11 in \cite{Brien06}
is based on the design with split-plots in a row--column design given
in Cochran and Cox~\cite{Cochran57}, Section 7.33.
Diagrams for the two randomizations are
given in Figure~\ref{fig:cc},
with leaf treatments named as viruses for clarity,
soil treatments designated as different soils for brevity,
and Altitude substituted
for Layer so that no two factors begin with the same letter.
Two diagrams are needed, because the assumed structure on leaves
changes between the randomizations, as shown in the two right-hand
panels.
At first sight, this experiment seems to involve
unrandomized-inclusive randomizations, because soils
are randomized to leaves in the first randomization, and then viruses are
randomized to leaves, taking into account the location of the soils.
However, the change in the assumed structure on the
leaves between the two randomizations makes them incoherent rather
than unrandomized-inclusive.

Table \ref{tab:cc} shows an attempt to build up a decomposition table
for this design. The first two columns follow directly from the
randomization in the top half of Figure~\ref{fig:cc}. The third
column corresponds to the leaves tier in the bottom half of
Figure~\ref{fig:cc}. When we use it to refine the decomposition given
by the first two tiers, we find that the Soils source occurs in two
tiers. Although this can happen in special circumstances like those
in Example \ref{eg:RotateGraze}, this is already a signal that
something may be wrong. We also find that the nesting, in this tier,
of Benches within Soils gives a source $\mathrm{B}[\mathrm{S}]$
with $8$ degrees of freedom. This is the sum of the previous sources
Benches and Residual in $\mathrm{Plants}[\mathrm{B}]$; these two parts
are denoted $\mathrm{B}[\mathrm{S}]_{\mathrm{B}}$ and
$\mathrm{B}[\mathrm{S}]_\vdash$ in Table \ref{tab:cc}.
Similarly, the nesting of Altitudes within Soils, in this tier, gives
sources $\mathrm{A}[\mathrm{S}]$ and
$\mathrm{A}\,\#\,\mathrm{B}[\mathrm{S}]$ which are each the sum of two
previous sources.

%
\begin{table}[b]
\caption{Proposed design for Example \protect\ref{eg:cc}
(columns denote plants; $s_0$--$s_3$ are different soils; $0$--$2$
denote viruses)}\label{tab:SPLwLS}
\begin{tabular*}{\tablewidth}{@{\extracolsep{\fill}}lcccccccccccc@{}}
\hline
& \multicolumn{4}{c}{\textbf{Bench I}} &
\multicolumn{4}{c}{\textbf{Bench II}}
& \multicolumn{4}{c@{}}{\textbf{Bench III}}\\[-4pt]
& \multicolumn{4}{c}{\hrulefill} &
\multicolumn{4}{c}{\hrulefill}
& \multicolumn{4}{c@{}}{\hrulefill}\\
\multicolumn{1}{@{}l}{\textbf{Soils}}
& \multicolumn{1}{c}{$s_0$} & \multicolumn{1}{c}{$s_1$} &
\multicolumn{1}{c}{$s_2$} & \multicolumn{1}{c}{$s_3$} &
\multicolumn{1}{c}{$s_0$} &
\multicolumn{1}{c}{$s_1$} & \multicolumn{1}{c}{$s_2$} &
\multicolumn{1}{c}{$s_3$} &
\multicolumn{1}{c}{$s_0$} & \multicolumn{1}{c}{$s_1$} &
\multicolumn{1}{c}{$s_2$} & \multicolumn{1}{c}{$s_3$}\\
\textbf{Altitude}\\
\hline
Top & 0 & 0 & 0 & 0 & 1 & 1 & 1 & 1 & 2 & 2 & 2 & 2\\
Middle & 2 & 2 & 1 & 1 & 0 & 0
& 2 & 2 & 1 & 1 & 0 & 0\\
Bottom & 1 & 1 & 2 & 2 & 2 & 2 & 0 & 0 & 0 & 0 & 1 & 1\\
\hline
\end{tabular*}
\end{table}

The real difficulties come when we try to incorporate the column for
the viruses tier, because the location of the Viruses source depends
on the outcome of the randomizations. For the outcome given in
\cite{Cochran57}, Section 7.33, and \cite{Brien06}, Example 11, the
Viruses source
does not have first-order balance in relation to either
$\mathrm{Altitudes}\,\#\,\mathrm{Benches}$ or
$\mathrm{Altitudes}\,\#\,\mathrm{Plants}[\mathrm{Benches}]$.
The interaction $\mathrm{Viruses}\,\#\,\mathrm{Soils}$
has the same problem.

If $\mathrm{Altitudes}\,\#\,\mathrm{Benches}$ is merged with
$\mathrm{Altitudes}\,\#\,\mathrm{Plants}[\mathrm{Benches}]$
in the decomposition table, then the analysis
is orthogonal and is equivalent to that given in
\cite{Cochran57}. However, this does not allow for consistent
Altitude differences across Plants, so it removes six spurious
degrees of freedom from what Cochran and Cox call ``Error (b)'' in
\cite{Cochran57}, page 310. The problem is that the design for
the Viruses does not respect the factor relationships
established in applying the Soils.
As Yates showed in \cite{Yates35}, if the randomization respects
Benches and Altitudes then a randomization-based model must include
their interaction.

What is needed is a design for a two-tiered experiment in which the
twelve treatments (combinations of levels of Soils and Viruses) are
randomized to leaves$_1$ in such a way that there is a refinement of
the natural decomposition of the treatments space which is
structure balanced in relation to Altitudes$\,\#\,$Benches.
For example, one might choose the systematic design in
Table \ref{tab:SPLwLS} and then randomize benches,
altitudes, and plants within benches.
In this design the twelve treatments are arranged in a $(3 \times
3)/4$ semi-Latin square constructed from a pair of mutually
orthogonal Latin squares of order $3$. The Viruses are arranged
according to one square for soils $s_0$ and $s_1$, and according to
the other square for $s_2$ and $s_3$.
Theorem 5.4 in \cite{Bailey92} shows that this design is the most
efficient with respect to $\mathrm{Altitudes}\,\#\,\mathrm{Benches}$.
Let $\mathrm{S}_1$ be a pseudofactor for Soils
whose two levels distinguish between the first two and the last two
levels of Soils. The design is structure balanced:
Viruses and $\mathrm{Viruses}\,\#\,\mathrm{S}_1$ have efficiency
factor $1/2$ in $\mathrm{Altitudes}\,\#\,\mathrm{Benches}$, while the
rest of the interaction $\mathrm{Viruses}\,\#\,\mathrm{Soils}$
is orthogonal to $\mathrm{Altitudes}\,\#\,\mathrm{Benches}$.
It has the advantage of having 10 degrees of freedom for the Residual for
$\mathrm{Altitudes}\,\#\,\mathrm{Plants}[\mathrm{Benches}]$,
two more than for the Cochran and Cox \cite{Cochran57} design.

Thus construction of the decomposition table when designing the
experiment can help to detect problems with a proposed design. In this
case, it helped to draw attention to the incoherent randomizations, to
highlight the associated problems and to give insight into how they
might be redressed.
\end{egg}

\subsection{Other structures}
\label{s:otherstruct}

All the examples in \cite{Brien09b} and this paper are poset block structures,
being defined by some factors and their nesting relationships,
as explained in \cite{Bailey04a,Brien09b}.
More generally, a structure may be a Tjur structure that is defined by
a family
of mutually orthogonal partitions or generalized factors (see \cite
{Tjur84} or
\cite{Bailey96}). Again, the generalized factors are summarized
in the Hasse diagram that depicts their marginality relations.
There is one projector $\mathbf{P}$ for each generalized factor $F$,
obtained from the Hasse diagram just as in Section 3
in \cite{Brien09b}, so that the effect of $\mathbf{P}$ on any vector
is still achieved by a straightforward sequence of averaging operations and
subtractions.
It is possible for some of these projectors to be zero.
Structures derived from tiers belong to this class.

Another common source of structure is an association scheme
\cite{Bose52,Bailey04a}: for example, the triangular scheme for all
unordered pairs from a set of parental types, which
is appropriate in a diallel experiment with no
self-crosses when the cross $(i,j)$ is regarded as the same as the
cross $(j,i)$. Then the matrices $\mathbf{P}$ are the minimal
idempotents of the association algebra \cite{Bose59}, and the corresponding
subspaces are its common eigenspaces \cite{Bailey04a}, Chapter 2.
The effect of $\mathbf{P}$ is a linear combination of the operations
of taking
sums over associate classes. In the case of the triangular association
scheme with $n$ parental types, the subspaces have dimensions $1$,
$n-1$ and $n(n-3)/2$; they correspond to the Mean, differences between
parental types and differences orthogonal to parental types, respectively.
The decomposition $\mathcal{R}_3$ in Example 5
in \cite{Brien09b} comes from
an association scheme with two associate classes.

The set of treatments in a rectangular lattice design exhibits yet
another kind of structure \cite{Bailey86}.
Although this structure derives neither from partitions nor
from an association scheme, the effect of each $\mathbf{P}$ is
achieved by averaging and subtracting.

The results here and in \cite{Brien09b}
apply to any structure that is an orthogonal decomposition
of the relevant vector space, so long as each structure can be
regarded as a decomposition of $V_\Omega$. For a Tjur structure
$\mathcal{Q}$ on a set $\Upsilon$ randomized to $\Omega$,
condition (4.1) in \cite{Brien09b}
must hold in order for $\mathcal{Q}$ to be
regarded as an orthogonal decomposition of $V_\Omega$. For structures
not defined by partitions, it seems that we need $\mathbf{Q}_i
\mathbf{X'X} \mathbf{Q}_j$ to be zero whenever $\mathbf{Q}_i \ne
\mathbf{Q}_j$,
where $\mathbf{X}$ is the $\Omega\times\Upsilon$ design matrix.
For an association scheme, this implies that the design
must be equireplicate. The analogue of
Theorem 5.1(a)
in \cite{Brien09b} for association schemes is given in
\cite{Bailey04a}, Section 7.7.

We admit that there are relevant experimental structures, such as
neighbour relations in a field or increasing quantities of dose, that
are not adequately described by an orthogonal decomposition of the
space. Nonetheless, a theory which covers designed experiments where
all the structures are orthogonal decompositions has wide
applicability, and we limit ourselves to such structures
here and in \cite{Brien09b}.

\subsection{Multiphase experiments}
Multiphase experiments are one of the commoner types of multitiered experiment.
As outlined in \cite{Brien06}, Section 8.1,
two-phase experiments may
involve almost any of the different types of multiple randomizations
and, as is evident from Section \ref{s:summary2}, these differ in their
assumptions.

If treatments are introduced only in the first phase,
then the randomizations form a chain, as in \cite{Brien09b}.
In \cite{Wood88}, Wood, Williams and Speed consider a
class of such two-phase designs for which
$\mathcal{R}$ is orthogonal in relation to
the natural structure $\mathcal{Q}_1$ on the middle tier,
and there is a refinement $\mathcal{Q}_2$ of $\mathcal{Q}_1$
such that $\mathcal{Q}_2 \vartriangleright\mathcal{R}$
is structure balanced in relation to $\mathcal{P}$.
The results there are less general than ours. First, the
assumptions for the second phase are in the nature of those for
randomized-inclusive randomizations only.
Second, the designs are restricted to those for
which the design for the first phase is orthogonal.

If treatments are introduced after the first phase, as in
Example \ref{eg:corn}, then some form of two-to-one randomization
is needed. Similarly,
Brien and\break Dem\'{e}trio \cite{Brien09} describe a three-phase experiment
involving composed and coincident randomizations.

\subsection{Further work}
While obtaining mixed model analyses of multitiered experiments
has been described in \cite{Brien06}, Section 7, and \cite{Brien09},
it remains to establish their randomization analysis.
The effects of intertier interactions on the analysis need to be investigated.
We would like to establish conditions under which
closed-form expressions are available for the Residual or Restricted Maximum
Likelihood (REML) estimates of the variance components \cite{Patterson71}
and Estimated
Generalized Least Squares (EGLS) estimates of the fixed effects.
Also required
is a derivation of the extended algorithm described in
\cite{Brien99} for obtaining the ANOVA for a multitiered
experiment.

Furthermore,
we have provided the basis for assessing a particular design for a
multitiered experiment, yet general principles for designing
them are still needed.

\printaddresses

\end{document}